\documentclass[journal]{IEEEtran}
\ifCLASSINFOpdf
   \usepackage[pdftex]{graphicx}
\else
\fi
\usepackage[cmex10]{amsmath}
\usepackage{amsfonts,amssymb}
\begin{document}
%
\title{Complete Real Time Solution of the General Nonlinear Filtering Problem without Memory}

%
%
%

\author{Xue~Luo~
        and~Stephen~S.-T.~Yau,~\IEEEmembership{Fellow,~IEEE}\\[9pt]
	{\large Dedicated to Professor Yu-Chi Ho and Professor Sanjoy Mitter on the occasion of their 80th birthday}%
\thanks{Manuscript received March 26, 2012; revised November 9, 2012 and March 9, 2013. This work is supported by the National Nature Science Foundation of China
(Grant No. 31271408) and the start-up fund from Tsinghua University.}
\thanks{X. Luo is with the department of
        Mathematics, Statistics and Computer Science,
        University of Illinois at Chicago, Science and Engineering Offices (M/C 249),
        851 S. Morgan Street, Chicago, IL 60607-7045 USA. (e-mail: xluo6@uic.edu).}
\thanks{S. S.-T. Yau is with the department of mathematical sciences, Tsinghua University, Beijing, 100084, P.R.China. (e-mail: yau@uic.edu).}
}

%
%

\markboth{IEEE TRANSACTIONS ON AUTOMATIC CONTROL,~Vol.~XX, No.~X, XXXX~2013}%
{Shell \MakeLowercase{\textit{et al.}}: Bare Demo of IEEEtran.cls for Journals}
%



\maketitle

\newtheorem{definition}{Definition}
\renewcommand{\thedefinition}{\arabic{section}.\arabic{definition}}
\newtheorem{proposition}{Proposition}
\renewcommand{\theproposition}{\arabic{section}.\arabic{proposition}}
\newtheorem{theorem}{Theorem}
\renewcommand{\thetheorem}{\arabic{section}.\arabic{theorem}}
\newtheorem{lemma}{Lemma}
\renewcommand{\thelemma}{\arabic{section}.\arabic{lemma}}
\newtheorem{corollary}{Corollary}
\renewcommand{\thecorollary}{\arabic{section}.\arabic{corollary}}
\newtheorem{remark}{Remark}
\renewcommand{\theremark}{\arabic{section}.\arabic{remark}}
\renewcommand{\theequation}{\arabic{section}.\arabic{equation}}
\newcommand{\inlineeq}{\hfill\refstepcounter{equation}{(\theequation)}}

\begin{abstract}
It is well known that the nonlinear filtering problem has important
applications in both military and civil industries. The central
problem of nonlinear filtering is to solve the
Duncan-Mortensen-Zakai (DMZ) equation in real time and in a
memoryless manner. In this paper, we shall extend the algorithm developed previously by S.-T. Yau and the second author to the most general setting of nonlinear filterings, where the explicit time-dependence is in the drift term, observation term, and the variance of the noises could be a matrix of functions of both time and the states. To preserve the off-line virture of the algorithm, necessary modifications are illustrated clearly. Moreover, it is shown rigorously that the approximated solution obtained by the algorithm converges to the real solution in the $L^1$ sense. And the precise error has been estimated. Finally, the numerical simulation support the feasibility and efficiency of our algorithm.
\end{abstract}


%
\IEEEpeerreviewmaketitle

\section{Introduction}

\setcounter{equation}{0}
%
%
%
%
\IEEEPARstart{T}{racing} back to 1960s, two most influential mathematics papers \cite{K}, \cite{KB} have been published in ASME's Journal of Basic Engineering. These are so-called Kalman filter (KF) and Kalman-Bucy filter. They addressed a significant question: How does one get accurate estimate from noisy data? The applications of KF are endless, from seismology to bioengineering to econometrics. The KF surpasses the other filtering in, at least, the following two aspects:  
\begin{itemize}
	\item The KF uses each new obervation to update a probability distribution for the state of the system without refering back to any earlier observations. This is so-called ``memoryless" or ``without memory".
	\item The KF makes the decisions of the state on the spot, while the observation data keep coming in. This property is called ``real time" application.
\end{itemize} 

Despite its success in many real applications, the limitations on the nonlinearity and Gaussian assumption of the initial probability density of the KF push the mathematicians and scientists to seek the optimal nonlinear filtering. One direction is to modify KF to adapt the nonlinearities. The researchers developed extended Kalman filter (EKF), unscented Kalman filter, ensemble Kalman filter, etc., which can handle weak nonlinearities (that is almost linear). But for serious nonlinearities, they may completely fail. The failure of EKF is shown in our numerical experiment, see the 1D cubic sensor in section VII.A.

Another direction, and also the most popular method nowadays, is the particle filter (PF), refer to as \cite{AMGC}, \cite{BC} and reference therein. It is developed from sequential Monte Carlo method. On the one hand, the PF is appliable to nonlinear, non-Gaussian state update and observation equations. As the number of particles goes to infinity, the PF estimate becomes asymptotically optimal. On the other hand, it is hard to be implemented as a real time application, due to its essence of Monte Carlo simulation. 

Besides the widely used two methods above, the partial differential equaitons (PDE) methods are introduced to the nonlinear filtering in 1960s. These methods are based on the fact that the unnormalized conditional density of the states is the solution of Duncan-Mortensen-Zakai (DMZ) equation, refer to as  \cite{Du}, \cite{Mo} and \cite{Za}. The classical PDE methods could be applied to this stochastic PDE to obtain an approximation to the density. Yet, the main drawback of PDE methods are the intensive computation. It is almost impossible to achieve the ``real time" performance. To overcome this shortcoming, the splitting-up algorithm is introduced to move the heavy computation off-line. It is like the Trotter product formula from semigroup theory. This operator splitting algorithm is proposed for the DMZ equation by Bensoussan, Glowinski, and Rascanu \cite{BGR}. More research articles follow this direction are \cite{GK}, \cite{Na} and \cite{It}, etc. Unfortunately, it is pointed out in \cite{BGR} that the soundness of this algorithm is verified only to the filtering with bounded drift term and observation term ( i.e., $f$ and $h$  in (\ref{Ito SDE})). Essentially with the similar idea, Yau and Yau \cite{YY} developed a novel algorithm to the ``pathwise-robust" DMZ equation (see (\ref{robust DMZ eqn expansion})), where the boundedness conditions are weakened by some mild growth conditions on $f$ and $h$. The two nice properties of the KF have also been kept in this algorithm: ``without memory" and ``real time". But their algorithm has only been rigorously proved in theory, when the drift term, the observation term ($f$, $h$ in \eqref{Ito SDE}) are not explicitly dependent on time, the variance of the noises ($G$ in \eqref{Ito SDE}) is the identity matrix, and the noises are standard Brownian motion processes ($S=I_{r\times r}$, $Q=I_{m\times m}$ in \eqref{Ito SDE}) .

In this paper, we shall extend the algorithm in \cite{YY} to the most general settings of nonlinear filtering problems, in the sense that the drift term, the observation term could explicitly depend on time, the variance of the noises $S$, $Q$ are time-dependent, and $G$ could be a matrix of functions of both time and the states. We shall validate our algorithm under very mild growth conditions on $f$, $h$ and $G$, see \eqref{condition c1}, \eqref{condition c3} and \eqref{N<c}. These are essentially time-dependent analogue of those in \cite{YY}. First of all, this extension is absolutely necessary. Many real applications have explicit time-dependence in their models, say the target orientation angles estimation from target position/velocity in constant turn model, where the angular velocities are piecewise constant functions in time \cite{R}. Second, this extension is nontrivial from the mathematical point of view. More trickier analysis of PDE is required. For instance, we need to take care of the more general elliptic operator $D_w^2$, see (2.5), rather than the Laplacian.

This paper is organized in the following. The detailed formulation of our algorithm is described in section II; In section III, we state our main theorems which validate our algorithm in theory. Notations and prelimilary are in section IV. Section V is devoted to the proofs of the main theorems. The lower bound of the density function is investigated in section VI. Numerical simulations are included in section VII. Finally, we arrive the conclusion. The appendices is consisted of the proof of the well-posedness theorem and the proof of an interesting property of the density function.

\section{Model and Algorithm}

\setcounter{equation}{0}

The model we are considering is the signal observation model with explicit time-dependence in the drift term, observation term and the variance of the noises:
\begin{equation}\label{Ito SDE}
   \left\{ \begin{aligned}
        dx_t &= f(x_t,t)dt+G(x_t,t)dv_t,\\
        dy_t &= h(x_t,t)dt+dw_t,
\end{aligned} \right.
\end{equation}
where $x_t$  and $f$ are $n$-vectors, $G$ is an $n\times r$ matrix,
and $v_t$ is an $r$-vector Brownian motion process with
$E[dv_tdv_t^{T}]=Q(t)dt$, $y_t$ and $h$ are $m$-vectors and $w_t$ is
an $m$-vector Brownian motion process with $E[dw_tdw_t^{T}]=S(t)dt$
and $S(t)>0$. We refer to $x_t$ as the state of the system at time
$t$ with some initial state $x_0$ (not necessarily obeying Gaussian distribution) and $y_t$ as the
observation at time $t$ with $y_0=0$. We assume that
$\{v_t,t\geq0\}$, $\{w_t,t\geq0\}$ and $x_0$ are independent. For the sake of convenience, let us call this system is the ``time-varying" case, while in \cite{YY} the ``time-invariant" case is studied.

Throughout this paper, we assume that $f$, $h$ and $G$ are $C^2$ in space and $C^1$ in time. Some growth conditions on $f$ and $h$ are expected to guaratee the existence and uniqueness of the ``pathwise-robust" DMZ equation.

The unnormalized density function $\sigma(x,t)$ of $x_t$ conditioned
on the observation history $Y_t=\{y_s:0\leq s\leq t\}$ satisfies the DMZ equation (for  the detailed formulation, see \cite{Du})
\begin{equation}\label{DMZ eqn}
   \left\{ \begin{aligned}
        d\sigma(x,t) &=
        L\sigma(x,t)dt+\sigma(x,t)h^{T}(x,t)S^{-1}(t)dy_t\\
        \sigma(x,0) &=\sigma_0(x),
\end{aligned} \right.
\end{equation}
where $\sigma_0(x)$ is the probability density of the initial
state $x_0$, and
\begin{align}\label{L}
    L(\ast) \equiv \frac12\sum_{i,j=1}^n\frac{\partial^2}{\partial x_i\partial
    x_j}\left[\left(GQG^{T}\right)_{ij}\ast\right]-\sum_{i=1}^n\frac{\partial(f_i\ast)}{\partial
    x_i}.
\end{align}

In this paper, we don't solve the DMZ equation directly, due to the following two reasons. On the one hand, the DMZ equation (\ref{DMZ eqn}) is a stochastic partial differential equation due to the term $dy_t$. There is no easy way
to derive a recursive algorithm to solve this equation. On the other hand, in real
applications, one may be more interested in constructing robust
state estimators from each observation path, instead of having certain statistical data of thousands of repeated experiments. Here, the robustness means our state esitmator  is not sensitive to the observation path. This property is important, since in most of the real applications,  the observation arrives and is processed at discrete moments in time. The state estimator is expected to still perform well based on the linear interpolation of the discrete observations, instead of the real continuous observation path. For each ``given" observation, making an invertible
exponential transformation \cite{Ro}
\begin{align}\label{Rozovsky's transformation}
    \sigma(x,t) = \exp{[h^{T}(x,t)S^{-1}(t)y_t]}\rho(x,t),
\end{align}
the DMZ equation is transformed into a deterministic partial
differential equation (PDE) with stochastic coefficients, which we will
refer as the ``pathwise-robust" DMZ equation
\begin{equation}\label{robust DMZ eqn}
   \left\{ \begin{aligned}
        \frac{\partial\rho}{\partial t}(x,t)&+\frac{\partial}{\partial
        t}(h^{T}S^{-1})^{T}y_t\rho(x,t)\\
        =&\exp{(-h^{T}S^{-1}y_t)}\left[L-\frac12h^{T}S^{-1}h\right]\\
		&\cdot[\exp{(h^{T}S^{-1}y_t)\rho(x,t)}]\\
        \rho(x,0) =& \sigma_0(x).
\end{aligned} \right.
\end{equation}
Or equivalently,
\begin{equation}\label{robust DMZ eqn expansion}
   \left\{ \begin{aligned}
        \frac{\partial\rho}{\partial t}(x,t) &=
        \frac12D^2_w\rho(x,t)+F(x,t)\cdot\nabla\rho(x,t)+J(x,t)\rho(x,t)\\
        \rho(x,0) &= \sigma_0(x),
\end{aligned} \right.
\end{equation}
where
\begin{align}\label{Dw^2}
    D^2_w =& \sum_{i,j=1}^n(GQG^{T})_{ij}\frac{\partial^2}{\partial x_i\partial
    x_j},\\\label{F}\notag
    F(x,t) =& \left[\sum_{j=1}^n\frac{\partial}{\partial x_j}\left(GQG^T\right)_{ij}\right.\\
           &\left.+\sum_{j=1}^n(GQG^{T})_{ij}\frac{\partial K}{\partial
            x_j}-f_i\right]_{i=1}^n,\\\label{J}\notag
    J(x,t) &=
        -\frac{\partial}{\partial t}(h^{T}S^{-1})^{T}y_t
        +\frac12\sum_{i,j=1}^n\frac{\partial^2}{\partial x_i\partial
            x_j}\left(GQG^T\right)_{ij}\\\notag
		&+\sum_{i,j=1}^n\frac{\partial}{\partial
            x_i}\left(GQG^T\right)_{ij}\frac{\partial K}{\partial
            x_j}\\\notag
            &+\frac12\sum_{i,j=1}^n(GQG^{T})_{ij}\left[\frac{\partial^2 K}{\partial x_i\partial x_j}
            +\frac{\partial K}{\partial x_i}\frac{\partial K}{\partial
            x_j}\right]\\
        &-\sum_{i=1}^n\frac{\partial f_i}{\partial
        x_i}-\sum_{i=1}^nf_i\frac{\partial K}{\partial
        x_i}-\frac12(h^{T}S^{-1}h),
\end{align}
in which
\begin{align}\label{K}
    K(x,t)=h^T(x,t)S^{-1}(t)y_t.
\end{align}

The existence and uniqueness of the ``pathwise-robust" DMZ equation \eqref{robust DMZ eqn expansion} has been investigated by Pardoux \cite{Pa}, Fleming-Mitter \cite{FM}, Baras-Blankenship-Hopkins \cite{BBH} and Yau-Yau \cite{YY05}, \cite{YY}. The well-posedness is guaranteed, when the drift term $f\in C^1$ and the observation term $h\in C^2$ are bounded in \cite{Pa}. Fleming and Mitter treated the case where $f$ and $\nabla f$ are bounded. Baras, Blankenship and Hopkins obtained the well-posedness result on the ``pathwise-robust" DMZ equation with a class of unbounded coefficients only in one dimension. In \cite{YY05}, Yau and Yau established the well-posedness result under the condition that $f$, $h$ have at most linear growth. In the appendices of \cite{YY}, Yau and Yau obtained the existence and uniqueness results in the weighted Sobolev space, where $f$ and $h$ satisfy some mild growth condition. It is necessary to point out that there is a gap in their proof of existence (Theorem A.4). In this paper, we circumvent the gap by more delicate analysis to give a time-dependent analogous well-posedness result to the ``pathwise-robust" DMZ equation under some mild growth conditions on $f$ and $h$ in Theorem \ref{existence and uniqueness}. 

The exact solution to (\ref{robust DMZ eqn}) or (\ref{robust DMZ eqn
expansion}), generally speaking, doesn't have a closed form. So many mathematicians pay their effort on seeking an efficient
algorithm to construct a good approximation. In this paper, we will extend the algorithm in \cite{YY} to the ``time-varying" case (cf. \eqref{Ito SDE}). We will not only give
the theoretical proof of the soundness of our algorithm, but also illustrate a ``time-varying" numerical
simulation to support our results. The difficulties are
two folds: on one hand, the well-posedness of the ``time-invariant" robust DMZ equation, under some conditions, has been investigated by  \cite{FM}, \cite{Pa}, \cite{YY}, etc., while that in the ``time-varying" case hasn't been established
yet; on the other hand, the ``time-varying" case will lead to more involved computations and more delicate analysis. For instance, the Laplacian in ``time-invariant" case is replaced by a time-dependent elliptic operator $D_w^2$ in \eqref{Dw^2}. Furthermore, the two nice properties of KF, namely ``memoryless" and ``real time", are preserved in our algorithm.

Let us assume that we know the observation time sequence $0=\tau_0<\tau_1<\cdots<\tau_k=T$ apriorily. But the observation data $\{y_{\tau_i}\}$ at each sampling time $\tau_i$, $i=0,\cdots,k$ are unknown until the on-line experiment runs. We call the computation ``off-line", if it can be performed without any on-line experimental data (or say pre-computed); otherwise, it is called ``on-line" computations. One only concerns the computational complexity of the on-line computations, since this hinges the success of ``real time" application. 
  
Let us denote the observation time sequence as $\mathcal{P}_k=\{0=\tau_0< \tau_1<\cdots<\tau_k=T\}$. Let $\rho_i$ be the solution of the robust DMZ
equation with $y_t=y_{\tau_{i-1}}$ on the interval $\tau_{i-1}\leq
t\leq\tau_i$, $i=1,2,\cdots,k$
\begin{equation}\label{robust DMZ eqn freezed}
   \left\{ \begin{aligned}
        \frac{\partial\rho_i}{\partial t}(x,t)&+\frac{\partial}{\partial
        t}\left(h^{T}S^{-1}\right)^{T}y_{\tau_{i-1}}\rho_i(x,t)\\
        =&\exp{\left(-h^{T}S^{-1}y_{\tau_{i-1}}\right)}\left[L-\frac12h^{T}S^{-1}h\right]\\
             &\cdot\left[\exp{\left(h^{T}S^{-1}y_{\tau_{i-1}}\right)\rho_i(x,t)}\right]\\
        \rho_1(x,0) &= \sigma_0(x),\\
        \textup{or}\phantom{\rho_1(x,0)}&\\
        \rho_i(x,\tau_{i-1}) =& \rho_{i-1}(x,\tau_{i-1}),\quad\textup{for}\
        i=2,3,\cdots,k.
\end{aligned} \right.
\end{equation}
Define the norm of $\mathcal{P}_k$ by $|\mathcal{P}_k|=\sup_{1\leq i\leq k}(\tau_i-\tau_{i-1})$. Intuitively, as $|\mathcal{P}_k|\rightarrow0$, we have
$$\sum_{i=1}^k\chi_{[\tau_{i-1},\tau_i]}(t)\rho_i(x,t)\rightarrow\rho(x,t)$$
in some sense, for all $0\leq t\leq T$, where $\rho(x,t)$ is the exact
solution of (\ref{robust DMZ eqn}). That is to say, intuitively, the denser the sampling time sequence is, the more accurate the approximate solution should be obtained. Even though the intuition is shown rigorously to be true, it is impractical to solve \eqref{robust DMZ eqn freezed} in the ``real time" manner, since the ``on-line" data $\{y_{\tau_i}\}$, $i=1,\cdots,k$, are contained in the coefficients of \eqref{robust DMZ eqn freezed}. Therefore, we have to numerically solve the time-consuming PDE on-line, every time after the new observation data coming in. Yet, the proposition below helps to move the heavy computations off-line. This is the key ingredient of the algorithm in \cite{YY}, so is in ours.
\begin{proposition}
    For each $\tau_{i-1}\leq t<\tau_i$, $i=1,2,\cdots,k$, $\rho_i(x,t)$ satisfies \textup{(\ref{robust DMZ eqn freezed})} if and only if
    \begin{equation}\label{Rozovsky's reverse transformation}
        u_i(x,t) = \exp{\left[h^{T}(x,t)S^{-1}(t)y_{\tau_{i-1}}\right]}\rho_i(x,t),
    \end{equation}
satisfies the Kolmogorov forward equation (KFE)
\begin{equation}\label{Kolmogorov equation}
    \frac{\partial u_i}{\partial t}(x,t) =
    \left(L-\frac12h^{T}S^{-1}h\right)u_i(x,t),
\end{equation}
where $L$ is defined in \textup{(\ref{L})}.
\end{proposition}

It is clear that \eqref{Kolmogorov equation} is independent of the observation path $\{y_{\tau_i}\}_{i=0}^k$, and the transformation between $u_i$ and $\rho_i$ is one-to-one. It is also not hard to see that \eqref{Kolmogorov equation} could be numerically solved beforehand. Observe that the operator $\left(L-\frac12h^{T}S^{-1}h\right)$ is time-varying, unlike that studied in \cite{YY}. Let us denote it as $U(t)$ for short and emphasis its time-dependence. But this doesn't affect the ``off-line" virture of our algorithm. Under certain conditions, $\{U(t)\}_{t\in[0,T]}$ forms a family of strong elliptic operators. Furthermore, the operator $U(t):D(U(t))\subset L^2(\mathbb{R}^n)\rightarrow L^2(\mathbb{R}^n)$ is the infinitesimal generator of the two-parameter semigroup $\mathcal{U}(t,\tau)$, for $t\geq\tau$. In particular, with the observation time sequence known $\{\tau_i\}_{i=1}^k$, we obtain a sequence of two-parameter semigroup $\{\mathcal{U}(t,\tau_{i-1})\}_{i=1}^k$, for $\tau_{i-1}\leq t<\tau_i$. Let us take the initial conditions of KFE (\ref{Kolmogorov equation}) at $t=\tau_i$ as a set of complete orthonormal base in $L^2(\mathbb{R}^n)$, say $\{\phi_l(x)\}_{l=1}^\infty$. We pre-compute the solutions of \eqref{Kolmogorov equation} at time $t=\tau_{i+1}$, denoted as $\{\mathcal{U}(\tau_{i+1},\tau_i)\phi_l\}_{l=1}^\infty$. These data should be stored in preparation of the on-line computations. Compared with the ``time-invariant" case, the price to pay is that the ``time-varying" case requires more storage capacity, since $\{\mathcal{U}(\tau_{i+1},\tau_i)\phi_l\}_{l=1}^\infty$ differs from each $\tau_i$, $i=1,\cdots,k$, and all of them need to be stored. In general, the longer simulation time is the more storage it requires in the ``time-varying" case. While the storage of the data is independent of the simulation time in the ``time-invariant" case. Nevertheless, it won't affect the off-line virture of our algorithm. 

The on-line computation in our algorithm is consisted of two parts at each time step $\tau_{i-1}$, $i=1,\cdots,k$. 
\begin{itemize}
	\item Project the initial condition $u_i(x,\tau_{i-1})\in L^2(\mathbb{R}^n)$ at $t=\tau_{i-1}$ onto the base $\{\phi_l(x)\}_{l=1}^\infty$, i.e., $u_i(x,\tau_{i-1})=\sum_{l=1}^\infty\hat{u}_{i,l}\phi_l(x)$. Hence, the
solution to (\ref{Kolmogorov equation}) at $t=\tau_i$ can be expressed as
\begin{align}\label{proj to basis}\notag
u_i(x,\tau_i)=&\mathcal{U}(\tau_i,\tau_{i-1})u_i(x,\tau_{i-1})\\
		=&\sum_{l=1}^\infty\hat{u}_{i,l}\left[\mathcal{U}(\tau_i,\tau_{i-1})\phi_l(x)\right],
\end{align}
where $\{\mathcal{U}(\tau_i,\tau_{i-1})\phi_l(x)\}_{l=1}^\infty$ have already been computed off-line.
	\item Update the initial condition of \eqref{Kolmogorov equation} at $\tau_i$ with the new observation $y_{\tau_i}$. Let us specify the observation updates (the initial condition of
(\ref{Kolmogorov equation}) ) for each time step. For $0\leq
t\leq\tau_1$, the initial condition is $u_1(x,0)=\sigma_0(x)$. At
time $t=\tau_1$, when the observation $y_{\tau_1}$ is available, 
\begin{align*}
    u_2(x,\tau_1)&\overset{\eqref{Rozovsky's reverse transformation}}=\exp{[h^{T}(x,\tau_1)S^{-1}(\tau_1)y_{\tau_1}]}\rho_2(x,\tau_1)\\
                 &\overset{\eqref{Rozovsky's reverse transformation}, \eqref{robust DMZ eqn freezed}}=\exp{[h^{T}(x,\tau_1)S^{-1}(\tau_1)y_{\tau_1}]}u_1(x,\tau_1),
\end{align*}
with the fact $y_0=0$. Here, $u_1(x,\tau_1)=\sum_{l=1}^\infty\hat{u}_{1,l}\left[\mathcal{U}(\tau_1,0)\phi_l(x)\right]$, where $\{\hat{u}_{1,l}\}_{l=1}^\infty$ is computed in the previous step, and $\{\mathcal{U}(\tau_1,0)\phi_l(x)\}_{l=1}^\infty$ are prepared by off-line computations.  Hence, we obtain the initial condition $u_2(x,\tau_1)$ of \eqref{Kolmogorov equation} for the next time interval $\tau_1\leq
t\leq\tau_2$. Recursively, the initial condition of \eqref{Kolmogorov equation} for $\tau_{i-1}\leq t\leq\tau_i$ is
\begin{align}\label{IC}\notag
    u_i(x,\tau_{i-1})=&\exp{[h^{T}(x,\tau_{i-1})S^{-1}(\tau_{i-1})(y_{\tau_{i-1}}-y_{\tau_{i-2}})]}\\
			&\cdot u_{i-1}(x,\tau_{i-1}),
\end{align}
for $i=2,3,\cdots,k$, where $u_{i-1}(x,\tau_{i-1})=\sum_{l=1}^\infty\hat{u}_{i-2,l}\left[\mathcal{U}(\tau_{i-1},\tau_{i-2})\phi_l(x)\right]$.
\end{itemize}
The approximation of $\rho(x,t)$, denoted as
$\hat{\rho}(x,t)$, is obtained
\begin{align}\label{hatrho}
    \hat{\rho}(x,t)=\sum_{i=1}^k\chi_{[\tau_{i-1},\tau_i]}(t)\rho_i(x,t),
\end{align}
where $\rho_i(x,t)$ is obtained from $u_i(x,t)$ by (\ref{Rozovsky's reverse
transformation}). And $\sigma(x,t)$ could be recovered by (\ref{Rozovsky's transformation}). 

A natural question comes
to us:
\begin{quote}
    {\sl Is $\hat{\rho}(x,t)$, obtained by our algorithm, a good
approximation of the exact solution $\rho(x,t)$ to \eqref{robust DMZ eqn}, for $(x,t)\in\mathbb{R}^n\times[0,T]$,
as $|\mathcal{P}_k|\rightarrow0$? If it is, then in what sense?}
\end{quote}

\section{Statements of the main theorems}
\setcounter{equation}{0}

In this section, we shall state the main theorems in this paper, which validate our algorithm in theory. Notice that $u_i$ in \eqref{Kolmogorov equation} and $\rho_i$ in \eqref{robust DMZ eqn freezed} are one-to-one. Hence, we shall deal with $\rho_i$ in the sequel. 

We first show that the exact solution $\rho$ of the ``pathwise-robust" DMZ equation \eqref{robust DMZ eqn} is well approximated by $\rho_R$ as $R\rightarrow\infty$, for any $t\in[0,T]$, where $\rho_R$ is the solution to \eqref{robust DMZ eqn} restricted on $B_R$ (the ball centered at the origin with the radius $R$):
\begin{equation}\label{robust DMZ eqn on ball}
   \left\{ \begin{aligned}
        \frac{\partial\rho_R}{\partial t}(x,t) =&
            \frac12D^2_w\rho_R(x,t)+F(x,t)\nabla\rho_R(x,t)\\
		&+J(x,t)\rho_R(x,t)\\
        \rho_R(x,0) =& \sigma_{0,B_R}(x)\\
        \rho_R(x,t)=&0\qquad\textup{for}\ (x,t)\in\partial
        B_R\times[0,T],
\end{aligned} \right.
\end{equation}
where $D_w^2$, $F(x,t)$ and $J(x,t)$ are defined in
(\ref{Dw^2})-(\ref{J}) and $\sigma_{0,\Omega}$ is defined as
\begin{equation}\label{sigma 0,omega}
   \sigma_{0,\Omega}(x)=
   \left\{ \begin{aligned}
        &\sigma_0(x),\qquad x\in\Omega_\epsilon\\
	&\textup{smooth},\qquad x\in\Omega\setminus\Omega_\epsilon\\
        &0,\qquad x\in\mathbb{R}^n\setminus\Omega,
\end{aligned} \right.
\end{equation}
in which $\Omega_\epsilon=\{x\in\Omega:\
\textup{dist}(x,\partial\Omega)>\epsilon\}$. Next, it is left to show that $\rho_R$ is well approximated by the solution obtained by our algorithm restricted on $B_R$. In fact, on the time interval $[0,\tau]$, $0<\tau\leq T$. Let us denote the time partition $\mathcal{P}_k^\tau=\{0=\tau_0<\tau_1<\cdots<\tau_k=\tau\}$. $\rho_R(x,\tau)$ is well approximated by $\rho_{k,R}(x,\tau)$, as $k\rightarrow+\infty$, in the $L^1$ sense, where $\rho_{k,R}$ is the solution of \eqref{robust DMZ eqn freezed} restricted on $B_R$.

For the notational convenience, let us denote
\begin{align}\label{N}\notag
    N(x,t)\equiv&-\frac{\partial}{\partial
    t}\left(h^TS^{-1}\right)y_t-\frac12D_w^2K\\
		&+\frac12D_wK\cdot\nabla
    K-f\cdot\nabla K-\frac12\left(h^TS^{-1}h\right),
\end{align}
where
\begin{align}\label{Dw}
        D_w\ast=\left[\sum_{j=1}^n\left(GQG^T\right)_{ij}(x,t)\frac{\partial\ast}{\partial x_j}\right]_{i=1}^n,
    \end{align}
and $D_w^2$ and $K$ are defined in (\ref{Dw^2}) and (\ref{K}), respectively.

The error estimate between $\rho$ and $\rho_R$ is given by the following theorem.
\begin{theorem}\label{convergence of rho R to rho}
    For any $T>0$, let $\rho(x,t)$ be a solution of the ``pathwise-robust" DMZ
    equation (\ref{robust DMZ eqn expansion}) in
    $\mathbb{R}^n\times[0,T]$. Let $R\gg1$ and $\rho_R$ be the solution to \eqref{robust DMZ eqn on ball}. Assume the following conditions are
    satisfied, for all $(x,t)\in\mathbb{R}^n\times[0,T]$:
\begin{enumerate}
    \item $N(x,t)+\frac32n\left|\left|GQG^{T}\right|\right|_\infty+|f-
        D_wK|\leq C,$\inlineeq{}\label{condition c1}
    \item $e^{-\sqrt{1+|x|^2}}\left[14n\left|\left|GQG^{T}\right|\right|_\infty+4\left|f-D_wK\right|\right]\leq
            \tilde{C},$\inlineeq{}\label{condition c3}
\end{enumerate}
    where $N$, $D_w$ and $K$ are defined in \eqref{N}, \eqref{Dw} and \eqref{K}, respectively, and $C$, $\tilde{C}$ are constants possibly depending on $T$. Let
$v=\rho-\rho_R$, then $v\geq0$ for all $(x,t)\in B_R\times[0,T]$ and
\begin{align}\label{estimation2 of thm convergence of rho R to rho}
    \int_{B_{\frac
    R2}}v(x,T)\leq
    \bar{C}e^{-\frac9{16}R}\int_{\mathbb{R}^n}e^{\sqrt{1+|x|^2}}\sigma_0(x),
\end{align}
where $\bar{C}$ is some constant, which may depend on $T$.
\end{theorem}

The next theorem tells us that $\rho_R$ is well approximated by the solution obtained by our algorithm restricted on $B_R$. More generally, $B_R$ is replaced by any bounded domain $\Omega\subset\mathbb{R}^n$ in the theorem.
\begin{theorem}\label{convergence of rho i to rho}
    Let $\Omega$ be a bounded domain in $\mathbb{R}^n$. Assume that
\begin{enumerate}
    \item $\left|N(x,t)\right|\leq C,$\inlineeq{}\label{N<c}
    \item
        There exists some $\alpha\in(0,1)$, such that
         \begin{align}\label{N-N<ct^alpha}
            \left|N(x,t)-N(x,t;\bar{t})\right|\leq \tilde{C}|t-\bar{t}|^\alpha,
        \end{align}
\end{enumerate}
    for all $(x,t)\in\Omega\times[0,T]$, $\bar{t}\in[0,T]$, where $N(x,t)$ is in \eqref{N}, and $N(x,t;\bar{t})$ denotes $N(x,t)$ with the observation $y_t=y_{\bar{t}}$. Let $\rho_\Omega(x,t)$ be the solution of (\ref{robust DMZ eqn expansion}) on $\Omega\times[0,T]$ with $0-$Dirichlet boundary condition:
    \begin{equation}
   \left\{ \begin{aligned}\label{rho omega}
        \frac{\partial\rho_\Omega}{\partial t}(x,t) =&
            \frac12D^2_w\rho_\Omega(x,t)+F(x,t)\cdot\nabla\rho_\Omega(x,t)\\
		&+J(x,t)\rho_\Omega(x,t)\\
        \rho_\Omega(x,0) =& \sigma_{0,\Omega}(x)\\
        \rho_\Omega(x,t)|_{\partial\Omega} =& 0,
\end{aligned} \right.
\end{equation}
where $D_w^2$, $F(x,t)$ and $J(x,t)$ are defined in
(\ref{Dw^2})-(\ref{J}) and $\sigma_{0,\Omega}$ is defined in \eqref{sigma 0,omega}. For any $0\leq\tau\leq T$, let
$\mathcal{P}_k^\tau=\{0=\tau_0<\tau_1<\tau_2<\cdots<\tau_k=\tau\}$ be a partition of $[0,\tau]$, where $\tau_i=\frac{i\tau}{k}$. Let
$\rho_{i,\Omega}(x,t)$ be the approximate solution obtained by our algorithm restricted on $\Omega\times [\tau_{i-1},\tau_i]$. Or equivalently, $\rho_{i,\Omega}$ is the solution on $\Omega\times [\tau_{i-1},\tau_i]$ of the equation 
\begin{equation}\label{rho i omega}
   \left\{ \begin{aligned}
        \frac{\partial\rho_{i,\Omega}}{\partial t}(x,t)
            =&\frac12D^2_w\rho_{i,\Omega}(x,t)+F(x,t;\tau_{i-1})\cdot\nabla\rho_{i,\Omega}(x,t)\\
            &+J(x,t;\tau_{i-1})\rho_{i,\Omega}(x,t)\\
        \rho_{i,\Omega}(x,\tau_{i-1}) =& \rho_{i-1,\Omega}(x,\tau_{i-1})\\
        \rho_{i,\Omega}(x,t)|_{\partial\Omega} =& 0,
\end{aligned} \right.
\end{equation}
for $i=1,2,\cdots,k$, with the convention that
$\rho_{1,\Omega}(x,0)=\sigma_{0,\Omega}(x)$. Here, $F(x,t;\tau_{i-1})$, $J(x,t;\tau_{i-1})$ denote $F(x,t)$, $J(x,t)$ with the observation $y_t=y_{\tau_{i-1}}$, respectively. Then
$$\rho_\Omega(x,\tau)=\lim_{k\rightarrow\infty}\rho_{k,\Omega}(x,\tau),$$
in the $L^1$ sense in space and the following estimate holds:
\begin{align}\label{L^1}
    \int_\Omega|\rho_\Omega-\rho_{k,\Omega}|(x,\tau)\leq
    \bar{C}\frac1{k^\alpha},
\end{align}
where $\bar{C}$ is a generic constant, depending on $T$, $\int_\Omega\sigma_{0,\Omega}$. The right-hand side of (\ref{L^1}) tends to zero as
$k\rightarrow\infty$.
\end{theorem}

\section{Notations and preliminary}

\setcounter{equation}{0}

Throughout the paper, let $\mathbb{Q}_T=\mathbb{R}^n\times[0,T]$. Let $H^1(\mathbb{R}^n)$ be the Sobolev space, equipped the norm
\[
    ||u(x)||_1^2 = \int_{\mathbb{R}^n}(u^2+|\nabla_x u|^2)dx.
\]
And let $H^{1;1}(\mathbb{Q}_T)$ be the functional space of both $t$ and $x$, with the norm
\[
    ||v(x,t)||_{1;1}^2 = \int_{\mathbb{Q}_T} (v^2+|\nabla_x
    v|^2+|\partial_tv|^2) dxdt.
\]
The subspace of $H^{1;1}(\mathbb{Q}_T)$ consisting of functions $v(x,t)$ which have compact supports in $\mathbb{R}^n$ for any $t$ is denoted as $H^{1;1}_0(\mathbb{Q}_T)$.
\begin{definition}
	The function $u(x,t)$ in $H_0^{1;1}(\mathbb{Q}_T)$ is called a {\it weak solution} of the initial value problem
	\begin{align*}
	\left\{\begin{aligned}
		&\sum_{i,j=1}^n\frac{\partial}{\partial x_i}\left(A_{ij}(x,t)\frac{\partial u}{\partial x_j}\right)+\sum_{i=1}^nB_i(x,t)\frac{\partial u}{\partial x_i}\\
		&\phantom{aaaa}+C(x,t)u=\frac{\partial u}{\partial t},\\
		&u(x,0)=u_0(x)
	\end{aligned}\right.
	\end{align*}
if for any function $\Phi(x,t)\in H_0^{1;1}(\mathbb{Q}_T)$ the following relation holds:
	\begin{align*}
		\int\int_{\mathbb{Q}_T}&\left[\sum_{i,j=1}^nA_{ij}\frac{\partial u}{\partial x_i}\frac{\partial\Phi}{\partial x_j}\right.\\
&\left.-\left(\sum_{i=1}^nB_i\frac{\partial u}{\partial x_i}+Cu+\frac{\partial u}{\partial t}\right)\Phi\right]dxdt=0
	\end{align*}
and $u(x,0)=u_0(x)$.
\end{definition}

We assume that the following conditions hold throughout the paper:
\begin{enumerate}
	\item The operator $L$ defined in (\ref{L}) is a strong elliptic
operator and it is bounded from above on $\mathbb{Q}_T$. That is, there
exists a constant $\lambda>0$ such that
\[
    \lambda|\xi|^2\leq\sum_{i,j=1}^n(GQG^T)_{ij}\xi_i\xi_j,
\]
for any $(x,t)\in \mathbb{Q}_T$, for any
$\xi=(\xi_1,\xi_2,\cdots,\xi_n)\in\mathbb{R}^n$. And
\[
    ||GQG^T||_\infty =
    \sup_{(x,t)\in\mathbb{Q}_T}|GQG^T|_\infty<\infty,
\]
where $|\cdot|_\infty$ is the sup-norm of the matrix.
	\item The initial density function $\sigma_0(x)\in
H^1(\mathbb{R}^n)$ decays fast enough. To be more specific, we
require that
\[
    \int_{\mathbb{R}^n}e^{\sqrt{1+|x|^2}}\sigma_0(x)<\infty.
\]
\end{enumerate}

\subsection{Existence and uniqueness of the non-negative weak solution}

For the sake of completeness, we state the existence and uniqueness of the non-negative weak solution to (\ref{robust DMZ eqn expansion}) on $\mathbb{Q}_T$ below. 
\begin{theorem}\label{existence and uniqueness}
    Under the conditions 1)-2) and the conditions (\ref{d/dt (GQG^T)})-(\ref{A condition 4}) in Theorem \ref{a priori estimate},
    the ``pathwise-robust" DMZ equation (\ref{robust DMZ eqn
    expansion}) on $\mathbb{Q}_T$ with the initial value $\sigma_0\in
    H^1(\mathbb{R}^n)$ admits a non-negative weak solution $\rho\in H^{1;1}(\mathbb{Q}_T)$. Assume further that the conditions \eqref{unique-cond1}-\eqref{alpha} in Theorem \ref{uniqueness} are satisfied, then the weak solution $\rho$ on $\mathbb{Q}_T$ is unique.
\end{theorem}

To avoid the distraction from our main theorems, we leave the detailed proof of this theorem in Appendix A.

\subsection{Technical lemma}

In the proofs of our main theorems, we will repeatedly adopt the
following lemma with suitably chosen test functions. Let us state the lemma and sketch the proof here.
\begin{lemma}\label{IBP}
    Assume that $\rho_\Omega$ satisfies the ``pathwise-robust" DMZ equation (\ref{robust DMZ eqn expansion}) on
    some bounded domain $\Omega\in\mathbb{R}^n$, for $0\leq t\leq T$. 
    Then, for any test function $\psi(x)\in C^\infty(\Omega)$, we have
    \begin{align}\label{lemma eq0}\notag
        \frac d{dt}\int_\Omega\psi\rho_\Omega
            =&\frac12\int_\Omega D_w^2\psi\rho_\Omega+\int_\Omega
                (f-D_wK)\cdot\nabla\psi\rho_\Omega\\\notag
	&+\int_\Omega\psi\rho_\Omega N+\frac12\int_{\partial\Omega}\psi\left(D_w\rho_\Omega\cdot\nu\right)\\\notag
       &-\frac12\int_{\partial\Omega}\rho_\Omega\left(D_w\psi\cdot\nu\right)\\\notag	&+\frac12\int_{\partial\Omega}\psi\rho_\Omega\sum_{i,j=1}^n\frac{\partial}{\partial
                x_i}\left(GQG^T\right)_{ij}\nu_j\\
             &+\int_{\partial\Omega}\psi\rho_\Omega\left(D_wK\cdot\nu\right)
                -\int_{\partial\Omega}\psi\rho_\Omega(f\cdot\nu),
    \end{align}
    where $\nu=(\nu_1,\nu_2,\cdots,\nu_n)$ is the exterior normal
    vector of $\Omega$, and $D_w^2$, $K$ and $D_w$ are defined in \eqref{Dw^2}, \eqref{K} and \eqref{Dw}, respectively.
\end{lemma}
\begin{IEEEproof}[Sketch of the proof]  Multiply $\psi(x)$ on both sides of (\ref{robust DMZ eqn expansion}) and integrate over
the domain $\Omega$, it yields
\begin{align}\label{lemma eq1}
    \frac d{dt}\int_\Omega\psi\rho_\Omega
        =&\int_\Omega\psi\left[\frac12D^2_w\rho_\Omega+F(x,t)\cdot\nabla\rho_\Omega+J(x,t)\rho_\Omega\right],
\end{align}
where $F(x,t)$ and $J(x,t)$ are defined in (\ref{F}) and (\ref{J}),
respectively. After applying integration by parts to the first two terms on
the right-hand side of (\ref{lemma eq1}), (\ref{lemma eq0}) is obtained by written in
compact notations.
\end{IEEEproof}

\section{Proofs of the main theorems}

\subsection{Reduction to the bounded domain case}

\setcounter{equation}{0}

In this section we shall prove that the solution $\rho$ to the
``pathwise-robust" DMZ equation (\ref{robust DMZ eqn expansion}) in
$\mathbb{R}^n$ can be well approximated by the solution $\rho_R$ of
 (\ref{robust DMZ eqn on ball}) in a large ball $B_R$. Moreover, the
 error estimate with respect to the radius $R$ is given explicitly in the $L^1$
 sense. Let $C$, $\tilde{C}$ and $\hat{C}$ denote the generic constants, which may differ from
 line to line.

We first show an interesting proposition, which reflects how the
density function in the large ball changing with respect to time. It
is also an important ingredient of the error estimate in Theorem \ref{convergence of rho R to rho}.

\begin{proposition}\label{rho concentrate at origin}
    For any $T>0$, let $\rho_R(x,t)$ be a solution of the ``pathwise-robust" DMZ
    equation restricted on $B_R$ \eqref{robust DMZ eqn on ball}. Assume
    that condition (\ref{condition c1}) is satisfied. Then
    \begin{equation}\label{estimation1 of thm concentration at
    origin}
        \int_{B_R}e^{\sqrt{1+|x|^2}}\rho_R(x,t)\leq
        e^{Ct}\int_{\mathbb{R}^n}e^{\sqrt{1+|x|^2}}\sigma_0(x).
    \end{equation}
\end{proposition}
\begin{IEEEproof} Choose the test function in Lemma
\ref{IBP} $\psi=e^{\phi_1}$ , where ${\phi_1}\in C^\infty(B_R)$,
$B_R=\{x\in\mathbb{R}^n:|x|\leq R\}$. Let $\rho_R$ be the solution
of the ``pathwise-robust" DMZ equation (\ref{robust DMZ eqn on ball}) on the
ball $B_R$. By Lemma \ref{IBP}, we have
\begin{align}\label{A1}\notag
        \frac{d}{dt}\int_{B_R} e^{\phi_1} \rho_R
        =&\int_{B_R}e^{\phi_1}\rho_R\left[\frac12\left(D_w^2{\phi_1}+D_w{\phi_1}\cdot\nabla{\phi_1}\right)\right.\\\notag
		&\left.+(f-D_wK)\cdot\nabla{\phi_1}+N\right]\\
         &+\frac12\int_{\partial B_R}e^{\phi_1}(D_w\rho_R\cdot\nu).
\end{align}
All the boundary integrals in (\ref{lemma eq0}) vanish, except the first
term in \eqref{lemma eq0}, since $\rho_R|_{\partial\Omega}=0$. Moreover, recall
that $\rho_R\geq0$ in $B_R$ and vanishes on $\partial B_R$ implies
that $\frac{\partial\rho_R}{\partial\nu}|_{\partial B_R}\leq0$.
Hence, on $\partial B_R$,
\begin{align*}
    (D_w\rho_R\cdot\nu)
        =&\sum_{i=1}^n\left[\sum_{j=1}^n(GQG^T)_{ij}\frac{\partial\rho_R}{\partial
        r}\frac{\partial r}{\partial x_j}\right]\nu_i\\\notag
        =&\frac{\partial\rho_R}{\partial
        r}\left[\sum_{i,j=1}^n(GQG^T)_{ij}\frac{x_j}r\frac{x_i}r\right]
        \leq0,
\end{align*}
by the positive definite assumption of $(GQG^T)$. Thus, (\ref{A1})
can be reduced further
\begin{align}\label{A11}\notag
    \frac{d}{dt}\int_{B_R} e^{\phi_1} \rho_R
        \leq\int_{B_R}e^{\phi_1}\rho_R&\left[\frac12\left(D_w^2{\phi_1}+D_w{\phi_1}\cdot\nabla{\phi_1}\right)\right.\\
	&\left.+(f-D_wK)\cdot\nabla{\phi_1}+N\right].
\end{align}
Choose ${\phi_1}(x)=\sqrt{1+|x|^2}$ and estimate the terms
containing ${\phi_1}$ on the right-hand side of (\ref{A11}) one by
one:
\begin{align}\label{estimate of D_w^2phi_1}\notag
    D_w^2{\phi_1}=&\sum_{i=1}^n\left(GQG^{T}\right)_{ii}\frac1{\sqrt{1+|x|^2}}\\\notag
	&-\sum_{i,j=1}^n\left(GQG^{T}\right)_{ij}\frac{x_ix_j}{(1+|x|^2)^{\frac32}}\\\notag
        \leq&\left|\left|GQG^{T}\right|\right|_\infty\left[\frac n{\sqrt{1+|x|^2}}+\frac{n|x|^2}{(1+|x|^2)^{\frac32}}\right]\\
        \leq&2n\left|\left|GQG^{T}\right|\right|_\infty,\\\notag
    D_w{\phi_1}\cdot\nabla{\phi_1}\label{estimate of D_wphi_1 nabla phi_1}
        =&\sum_{i,j=1}^n\left(GQG^{T}\right)_{ij}\frac{x_ix_j}{1+|x|^2}\\
	\leq&\left|\left|GQG^{T}\right|\right|_\infty\frac{\sum_{i,j=1}^nx_ix_j}{1+|x|^2}
        \leq n\left|\left|GQG^{T}\right|\right|_\infty,
\end{align}
and
\begin{align}\label{estimate of (f-D_wK)nabla phi_1}\notag
    |(f-D_wK)\cdot\nabla\phi_1|\leq&|f-D_wK|\cdot\frac{|x|}{\sqrt{1+|x|^2}}\\
\leq&|f-D_wK|,
\end{align}
where $|\cdot|$ is the Euclidean norm. Substitute the estimate
(\ref{estimate of D_w^2phi_1})-(\ref{estimate of (f-D_wK)nabla
phi_1}) back into (\ref{A11}), we get
\begin{align*}
        \frac{d}{dt}&\int_{B_R} e^{\phi_1} \rho_R\\
        \leq&\int_{B_R}e^{\phi_1}\rho_R\left[\frac32n\left|\left|GQG^{T}\right|\right|_\infty+|f-D_wK|+N\right]\\
        \leq& C\int_{B_R}e^{\phi_1}\rho_R,
\end{align*}
by condition (\ref{condition c1}). Hence,
\begin{align*}
    \int_{B_R}e^{\phi_1}\rho_R(x,t)\leq&
    e^{Ct}\int_{B_R}e^{\phi_1}\rho_R(x,0)\leq
    e^{Ct}\int_{\mathbb{R}^n}e^{\phi_1}\rho(x,0)\\
	=&e^{Ct}\int_{\mathbb{R}^n}e^{\phi_1}\sigma_0(x),
\end{align*}
for $0\leq t\leq T$.
\end{IEEEproof}

We are ready to show Theorem \ref{convergence of rho R to rho}, i.e. the solution $\rho_R$ to (\ref{robust DMZ
eqn on ball}) on $B_R$ is a good approximation of $\rho$, the solution to \eqref{robust DMZ eqn expansion} in $\mathbb{R}^n$.

\begin{IEEEproof}[Proof of Theorem \ref{convergence of rho R to
rho}] By the maximum principle (cf. Theorem 1, \cite{Fr}), we have $v=\rho-\rho_R\geq0$ for
$(x,t)\in B_R\times[0,T]$, since $v|_{\partial B_R}\geq0$ for $0\leq
t\leq T$. Let us choose $\psi$ in Lemma \ref{IBP} as
$$\varrho(x)=e^{-{\phi_2}(x)}-e^{-R},$$
where ${\phi_2}$ is a radial symmetric function such that
${\phi_2}(x)|_{\partial B_R}=R$, $\nabla{\phi_2}|_{\partial B_R}=0$
and ${\phi_2}$ is increasing in $|x|$. Hence, $\varrho|_{\partial
B_R}=0$ and $\nabla\varrho|_{\partial B_R}=0$. Apply Lemma \ref{IBP}
to $v$, taking the place of $\rho_\Omega$, with the test function
$\psi=\varrho$, we have
\begin{align*}
    \frac{d}{dt}\int_{B_R}\varrho v
         =& \int_{B_R}v\left[\frac12D_w^2\varrho+\left(f-D_wK\right)\cdot\nabla
            \varrho+\varrho N\right]\\
         =&\int_{B_R}v\left\{\frac12e^{-{\phi_2}}\left(D_w{\phi_2}\cdot\nabla{\phi_2}-D_w^2{\phi_2}\right)\right.\\
            &\phantom{\int_{B_R}vaa}\left.-e^{-{\phi_2}}\left(f-D_wK\right)\cdot\nabla{\phi_2}+\varrho N\right\}\\
         =&\int_{B_R}v\varrho\left[-\frac12D_w^2{\phi_2}+\frac12D_w{\phi_2}\cdot\nabla{\phi_2}\right.\\
	&\phantom{\int_{B_R}vaa}\left.-\left(f-D_wK\right)\cdot\nabla{\phi_2}+N\right]\\
         &+e^{-R}\int_{B_R}e^{\sqrt{1+|x|^2}}v\left[e^{-\sqrt{1+|x|^2}}\right.\\
&\phantom{e^{-R}\int}\cdot\left(-\frac12D_w^2{\phi_2}+\frac12D_w{\phi_2}\cdot\nabla{\phi_2}\right.\\
            &\phantom{e^{-R}\int_{B_R}\cdot(}\left.\left.-\left(f-D_wK\right)\cdot\nabla{\phi_2}\right)\right]\\
         \triangleq&\int_{B_R}v\varrho\mathrm{I_1}+e^{-R}\int_{B_R}e^{\sqrt{1+|x|^2}}v\mathrm{I_2}.
\end{align*}
Let us choose $\phi_2(x)$ in $\varrho(x)$ to be
$\phi_2(x)=R\vartheta(\frac{|x|^2}{R^2})$, where
$\vartheta(x)=1-(1-x)^2$. It is easy to check that $\phi_2(x)$
satisfies all the conditions we mentioned before. Direct
computations yield, for any $x\in B_R$, $R>>1$,
\begin{align}\label{estimate of D_w^2phi_2}\notag
    \left|D_w^2{\phi_2}\right|
    =& \left|\sum_{i,j=1}^n\left(GQG^{T}\right)_{ij}\left(-\frac{8x_ix_j}{R^3}\right)\right.\\\notag
        &\left.+\sum_{i=1}^n\left(GQG^{T}\right)_{ii}\frac4R\left(1-\frac{|x|^2}{R^2}\right)\right|\\
    \leq&
    \left|\left|GQG^{T}\right|\right|_\infty\left(\frac{8n|x|^2}{R^3}+\frac{4n}R\right)
        \leq12n\left|\left|GQG^{T}\right|\right|_\infty,
\end{align}
\begin{align}\label{estimate of D_wphi_2 nabla
        phi_2}\notag
    \left|D_w{\phi_2}\cdot\nabla{\phi_2}\right|
    =&\left|\left(1-\frac{|x|^2}{R^2}\right)^2\sum_{i,j=1}^n\left(GQG^{T}\right)_{ij}\frac{4x_i}R\frac{4x_j}R\right|\\
    \leq&16n\left|\left|GQG^{T}\right|\right|_\infty,
\end{align}
and
\begin{align}\label{estimate of (f-D_wK)nabla phi_2}\notag
    \left|\left(f-D_wK\right)\cdot\nabla{\phi_2}\right|
    =& \left|\left(f-D_wK\right)\frac{4x}R\left(1-\frac{|x|^2}{R^2}\right)\right|\\
        \leq&4\left|f-D_wK\right|.
\end{align}
It follows that
\begin{align*}
    \sup_{B_R}|\mathrm{I_1}|
    \leq 14n\left|\left|GQG^{T}\right|\right|_\infty+4\left|f-D_wK\right|+N
    \leq C,
\end{align*}
by condition (\ref{condition c1}). Similarly,
\begin{align*}
    \sup_{B_R}|\mathrm{I_2}|
    \leq&\sup_{B_R}\left[e^{-\sqrt{1+|x|^2}}\left(14n\left|\left|GQG^{T}\right|\right|_\infty\right.\right.\\
        &\phantom{\sup_{B_R}[e^{-\sqrt{1+|x|^2}}(}
\left.\left.+4\left|f-D_wK\right|\right)\right]\\
    \leq& \tilde{C},
\end{align*}
by condition (\ref{condition c3}). In the view of Proposition
\ref{rho concentrate at origin}, one gets
\begin{align}\label{Thm 1.2 eqn 1}\notag
    \frac{d}{dt}\int_{B_R}\varrho v
        \leq& C\int_{B_R}\varrho v+e^{-R}\tilde{C}\int_{B_R}e^{\sqrt{1+|x|^2}}\rho\\
        \leq& C\int_{B_R}\varrho
            v+e^{-R+\hat{C}T}\tilde{C}\int_{\mathbb{R}^n}e^{\sqrt{1+|x|^2}}\sigma_0(x).
\end{align}
Multiply $e^{-Ct}$ on both sides of (\ref{Thm 1.2 eqn 1}) yields
\begin{align*}
    \frac d{dt}\left[e^{-Ct}\int_{B_R}\varrho v\right]
        \leq
        e^{-R+\hat{C}T-Ct}\tilde{C}\int_{\mathbb{R}^n}e^{\sqrt{1+|x|^2}}\sigma_0(x).
\end{align*}
Integrate from $0$ to $T$ and multiply $e^{CT}$ on both sides gives
us
\begin{align*}
    \int_{B_R}\varrho v(x,T)
        \leq&||v(x,0)||_\infty e^{CT}\int_{B_R}\varrho
        dx\\
	&+\frac{e^{CT}-1}Ce^{-R+\hat{C}T}\tilde{C}\int_{\mathbb{R}^n}e^{\sqrt{1+|x|^2}}\sigma_0(x),
\end{align*}
where $v(x,0)=\sigma_0-\sigma_{0,R}$. Recall that
$\varrho(x)=e^{-R\left[-\left(|x|^2/R^2-1\right)^2+1\right]}-e^{-R}$,
$|x|\leq R$, we arrive the following estimates:
\begin{align*}
    \int_{B_R}\varrho\leq\int_{B_R}\left(1-e^{-R}\right)\leq CR^n
\end{align*}
and
\begin{align*}
    \int_{B_R}\varrho v(x,T)
        \geq& \int_{B_{\frac
        R2}}\left(e^{-R\left[-\left(|x|^2/R^2-1\right)^2+1\right]}-e^{-R}\right)v(x,T)\\
        \geq&\frac12e^{-\frac7{16}R}\int_{B_{\frac R2}}v(x,T).
\end{align*}
It is easy to see that $||v(x,0)||_\infty\int_{B_R}\varrho\leq
C(n)\epsilon R^n$ is arbitrarily small, since $\epsilon$ is
independent of $R$. It yields that
\begin{align*}\tag{1.28}
\int_{B_{\frac R2}}v(x,T)\leq Ce^{-\frac9{16}R}
\int_{\mathbb{R}^n}e^{\sqrt{1+|x|^2}}\sigma_0(x),
\end{align*}
where $C$ is a generic constant, depending on $T$.
\end{IEEEproof}

By refining the proofs of Proposition \ref{rho concentrate at
origin} and Theorem \ref{convergence of rho R to rho}, we obtain an
interesting property of the density function $\rho(x,t)$. It asserts
that $\rho$ captures almost all the density in a large ball. And we
could give an precise estimate of the density outside the large
ball.

\begin{theorem}\label{coro}
Let $\rho(x,t)$ be a solution of the ``pathwise-robust" DMZ
    equation (\ref{robust DMZ eqn expansion}) in
    $\mathbb{Q}_T$. Assume that
    \begin{enumerate}
        \item Condition (\ref{condition c1}) is satisfied;
        \item A stronger version of condition (\ref{condition c3})
        is valid. To be more precise,
        \begin{align}\label{condition c}
            e^{-\frac12\sqrt{1+|x|^2}}\left[16n\left|\left|GQG^{T}\right|\right|_\infty+4\left|f-D_wK\right|\right]\leq C,
        \end{align}
     for all $(x,t)\in\mathbb{Q}_T$.
    \end{enumerate}
Then
    \begin{equation}\label{estimation2 of thm concentration at
    origin}
        \int_{|x|\geq R}\rho(x,T)
        \leq
        Ce^{-\frac12\sqrt{1+R^2}}\int_{\mathbb{R}^n}e^{\sqrt{1+|x|^2}}\sigma_0(x),
    \end{equation}
    where $C$ is a generic constant, which depends on $T$.
\end{theorem}

To avoid the distraction, we leave the detailed proof in Appendix B.

\subsection{$L^1$ convergence}

\setcounter{equation}{0}

In this section, we shall show that, for any $0<\tau\leq T$, with the partition $\mathcal{P}_k^\tau=\{0=\tau_0<\tau_1<\cdots<\tau_k=\tau\}$, the $L^1$ convergence of $\rho_{k,R}(x,\tau)$ to $\rho_R(x,\tau)$ holds, as $k\rightarrow+\infty$,  where $\rho_{k,R}$ is the solution of \eqref{rho i omega} obtained by our algorithm, and $\rho_R$ is the solution to \eqref{robust DMZ eqn on ball}. For the clarity, we state the technique lemma will be used in the proof of
Theorem \ref{convergence of rho i to rho} below.
\begin{lemma}\textup{(Lemma 4.1, \cite{YY})}\label{lemma 4.1
YauYau2}
    Let $\Omega$ be a bounded domain in $\mathbb{R}^n$ and let
    $v:\overline{\Omega}\times[0,T]\rightarrow\mathbb{R}$ be a $C^1$
    function. Assume that $v(x,t)=0$ for
    $(x,t)\in\partial\Omega\times[0,T]$. Let
    $\Omega_t^+=\{x\in\Omega:v(x,t)\geq0\}$. Then
    \[
        \frac{d}{dt}\int_{\Omega_t^+}v(x,t)=\int_{\Omega_t^+}\frac{\partial v}{\partial
        t}(x,t),
    \]
    for almost all $t\in[0,T]$.
\end{lemma}
\begin{IEEEproof}[ Proof of Theorem \ref{convergence of rho i to
rho}] For the notational convenience, we omit the subscript
$\Omega$ for $\rho_\Omega$ and $\rho_{i,\Omega}$ in this proof. Let
$\Omega_t^+=\{x\in\Omega:\rho(x,t)-\rho_i(x,t)\geq0\}$. Apply Lemma
\ref{IBP} to $(\rho-\rho_i)$ taking place of $\rho_\Omega$, with the test
function $\psi\equiv1$, we have
\begin{align}\label{estimate1 in thm C}\notag
    \frac{d}{dt}\int_{\Omega_t^+}(\rho-\rho_i)
    \leq&\int_{\Omega_t^+}(\rho-\rho_i)N(\cdot,t)\\
	&+\int_{\Omega_t^+}\rho_i[N(\cdot,t)-N(\cdot,t;\tau_{i-1})],
\end{align}
by Lemma \ref{lemma 4.1 YauYau2}. All the boundary integrals vanish,
except $\int_{\partial\Omega_t^+}D_w(\rho-\rho_i)\cdot\nu$, since
$(\rho-\rho_i)|_{\partial\Omega_t^+}=0$. Moreover,
$\int_{\partial\Omega_t^+}D_w(\rho-\rho_i)\cdot\nu\leq0$, due to the similar
argument for $\int_{\partial B_R}D_w\rho\cdot\nu\leq0$ in Proposition \ref{rho concentrate at origin}. Combine the conditions
(\ref{N<c}) and (\ref{N-N<ct^alpha}), (\ref{estimate1 in thm C}) can
be controlled by
\begin{align}\label{estimate2 in Thm C}
    \frac{d}{dt}\int_{\Omega_t^+}(\rho-\rho_i)
        \leq C\int_{\Omega_t^+}(\rho-\rho_i)+\tilde{C}(t-\tau_{i-1})^\alpha\int_\Omega\rho.
\end{align}
To estimate $\int_\Omega\rho$, we apply Lemma \ref{IBP} to $\rho$,
with the test function $\psi\equiv1$, we get
\begin{align*}
    \frac d{dt}\int_\Omega\rho
        \leq\int_\Omega\rho N
        \leq C\int_\Omega\rho,
\end{align*}
which implies
\begin{align}\label{int rho}
    \int_\Omega\rho
        \leq C\int_\Omega\sigma_{0,\Omega},
\end{align}
where $C$ is a generic constant, depending on $T$, for all $0\leq
t\leq T$. Thus,
\[
    \frac d{dt}\int_{\Omega_t^+}(\rho-\rho_i)
        \leq
        C\int_{\Omega_t^+}(\rho-\rho_i)+\tilde{C}(t-\tau_{i-1})^\alpha\int_\Omega\sigma_{0,\Omega}.
\]
Multiply $e^{-\tilde{C}(t-\tau_{i-1})}$ on both sides and integrate
from $\tau_{i-1}$ to $t$, we get
\begin{align*}
    \int_{\Omega_t^+}(\rho-\rho_i)(x,t)
        \leq& e^{\tilde{C}(t-\tau_{i-1})}\int_{\Omega_{\tau_{i-1}}^+}(\rho-\rho_i)(x,\tau_{i-1})\\
           &+C\frac{(t-\tau_{i-1})^{1+\alpha}}{1+\alpha}e^{\tilde{C}(t-\tau_{i-1})},
\end{align*}
where $C$ is a constant, which depends on $T$,
$\int_\Omega\sigma_{0,\Omega}$. Similarly, one can also get, for $\Omega_t^-=\{x\in\Omega:\,\rho(x,t)-\rho_i(x,t)<0\}$, that
\begin{align*}
    \int_{\Omega_t^-}(\rho_i-\rho)(x,t)
        \leq& e^{\tilde{C}(t-\tau_{i-1})}\int_{\Omega_{\tau_{i-1}}^-}(\rho_i-\rho)(x,\tau_{i-1})\\
            &+C\frac{(t-\tau_{i-1})^{1+\alpha}}{1+\alpha}e^{\tilde{C}(t-\tau_{i-1})}.
\end{align*}
Consequently, we have
\begin{align}\label{estimate of rho-rho_i}\notag
    &\int_\Omega|\rho-\rho_i|(x,t)\\\notag
        \leq&e^{\tilde{C}(t-\tau_{i-1})}\left[\int_\Omega|\rho-\rho_i|(x,\tau_{i-1})
            +C\frac{(t-\tau_{i-1})^{1+\alpha}}{1+\alpha}\right]\\
        \leq&e^{\tilde{C}(t-\tau_{i-1})}\left[\int_\Omega|\rho-\rho_{i-1}|(x,\tau_{i-1})
            +C\frac{(t-\tau_{i-1})^{1+\alpha}}{1+\alpha}\right],
\end{align}
since $\rho_i(x,\tau_{i-1})=\rho_{i-1}(x,\tau_{i-1})$, for
$i=1,2,\cdots,k$. Applying (\ref{estimate of rho-rho_i})
recursively, we obtain
\begin{align*}
    &\int_\Omega|\rho-\rho_k|(x,\tau_k)\\
        \leq&
        e^{\tilde{C}(\tau_k-\tau_{k-1})}\left[\int_\Omega|\rho-\rho_{k-1}|(x,\tau_{k-1})+C\frac{(\tau_k-\tau_{k-1})^{1+\alpha}}{1+\alpha}\right]\\
        \leq&e^{\tilde{C}T}\int_\Omega|\rho-\rho_0|(x,0)\\
            &+\frac{C}{1+\alpha}[(\tau_k-\tau_{k-1})^{1+\alpha}e^{\tilde{C}(\tau_k-\tau_{k-1})}\\
            &\phantom{+\frac{2c_6}{1+\alpha}[}+(\tau_{k-1}-\tau_{k-2})^{1+\alpha}e^{\tilde{C}(\tau_k-\tau_{k-2})}\\
	&\phantom{+\frac{2c_6}{1+\alpha}[}+\cdots+(\tau_1-\tau_0)^{1+\alpha}e^{\tilde{C}(\tau_k-\tau_0)}]\\
        =&\frac{C}{1+\alpha}\frac{T^{1+\alpha}}{k^{1+\alpha}}\left(e^{\tilde{C}\frac Tk}+e^{\tilde{C}\frac{2T}k}+\cdots+e^{\tilde{C}\frac{kT}k}\right)
            \leq\frac{C}{k^\alpha},
\end{align*}
where $C$ is a constant, which depends on $\alpha$, $T$ and
$\int_\Omega\sigma_{0,\Omega}$. It is clear that
$\int_\Omega|\rho-\rho_k|\rightarrow0$, as
$k\rightarrow\infty$.
\end{IEEEproof}

\section{Lower bound estimate of density function}

\setcounter{equation}{0}

It is well-known that solving the ``pathwise-robust" DMZ equation numerically is not easy because it is easily vanishing. We are also interested in whether the lower bound of the density function could be derived in the case where the drift term $f$ and the observation term $h$ are with at most the polynomial growth. The theorem below gives this lower bound:

\begin{theorem}\label{lower estimate}
    Let $\rho_R$ be the solution of (\ref{robust DMZ eqn on ball}),
    the ``pathwise-robust" DMZ equation on $B_R$. Assume that
\begin{enumerate}
    \item $f(x,t)$ and $h(x,t)$ have at most polynomial
        growth in $|x|$, for all $t\in[0,T]$;
    \item For any $0\leq t\leq T$, there exists positive
        integer $m$ and positive constants $C'$ and $C''$ independent
        of $R$ such that the following two conditions hold on
        $\mathbb{R}^n$:
        \begin{align}\label{condition c'}\notag
            \textup{(a)}\quad&\frac{|x|^{m-2}}2\left[nm(m-2)\left|\left|GQG^{T}\right|\right|_\infty\right.\\\notag
	&\phantom{\frac{|x|^{m-2}}2[}\left.+m\,\textup{Tr}\left(GQG^{T}\right)\right]\\
	&-m|x|^{m-2}(f-D_wK)\cdot x
            +N(x,t)\geq-C';\\\label{condition c''}\notag
            \textup{(b)}\quad&\left|n\left|\left|GQG^{T}\right|\right|_\infty\left(\frac12m^2|x|^{2m-2}\right.\right.\\\notag
		&\phantom{|n\left|\left|GQG^{T}\right|\right|_\infty(}\left.-m\left(\frac12m-1\right)|x|^{m-2}\right)\\\notag
	&-\frac12m\,\textup{Tr}\left(GQG^T\right)|x|^{m-2}\\\notag
	&\left.-m(f-D_wK)\cdot
                x|x|^{m-2}\right|\\
\leq&\frac12nm(m+1)\left|\left|GQG^{T}\right|\right|_\infty|x|^{2m-2}+C'',
        \end{align}
        where $\textup{Tr}(*)$ is the trace of $*$.
    \item Condition (\ref{N<c}) is satisfied.
\end{enumerate}
Then for any $R_0<R$,
    \begin{align*}
        &\int_{B_{R_0}}\zeta\rho_R(x,T)\\
        \geq&\frac{e^{(C-C')T-R_0^m}}{C'}\left(\frac12nm(m+1)\left|\left|GQG^T\right|\right|_\infty
    R_0^{2m-2}+C''\right)\\
	&\cdot\left(1-e^{C'T}\right)\int_{B_R}\sigma_{0,R}(x)
		+e^{-C'T}\int_{B_{R_0}}\zeta\sigma_{0,R}(x),
    \end{align*}
where $\zeta(x)=e^{-\xi(x)}-e^{-\xi(R_0)}$, $\xi(x)=|x|^m$.

    In particular, the solution $\rho$ of the ``pathwise-robust" DMZ equation (\ref{robust DMZ eqn}) on
    $\mathbb{R}^n$ has the estimate
    \[
        \int_{\mathbb{R}^n}e^{-|x|^m}\rho(x,T)\geq
        e^{-C'T}\int_{\mathbb{R}^n}e^{-|x|^m}\sigma_0(x).
    \]

\end{theorem}
\begin{IEEEproof} Apply Lemma \ref{IBP} to $\rho_R$ with
the test function $\psi$ to be $\zeta=e^{-\xi(x)}-e^{-\xi(R_0)}$,
where $\xi(x)$ is an increasing function in $|x|$, we have
\begin{align*}
    \frac{d}{dt}\int_{B_{R_0}}\zeta\rho_R
        =\int_{B_{R_0}}\rho_R\left[\frac12D_w^2\zeta+\left(f-D_wK\right)\cdot\nabla
            \zeta+\zeta N\right].
\end{align*}
All the boundary integrals vanish, since $\zeta|_{\partial
B_R}=\rho_R|_{\partial B_R}=0$. Direct computations yield that
\begin{align*}
    &\frac{d}{dt}\int_{B_{R_0}}\zeta\rho_R\\
    =&\int_{B_{R_0}}\rho_Re^{-\xi(R_0)}\\
	&\cdot\left\{\frac12\frac{\xi'^2(r)}{r^2}\sum_{i,j=1}^n\left(GQG^T\right)_{ij}x_ix_j
	 -\frac{\xi'(r)}r(f-D_wK)\cdot x\right.\\
            &\phantom{aa}           -\frac12\sum_{i,j=1}^n\left(GQG^T\right)_{ij}\left[\left(\xi''(r)-\frac{\xi'(r)}r\right)\frac{x_ix_j}{r^2}\right]\\
		&\phantom{aa}\left.
            -\frac12\textup{Tr}\left(GQG^T\right)\frac{\xi'(r)}r\right\}\\
            &+\int_{B_{R_0}}\zeta\rho_R\left[\frac12D_w\xi\cdot\nabla\xi-\frac12D_w^2\xi-(f-D_wK)\cdot\nabla\xi+N\right]\\
    \triangleq&\mathrm{I_3}+\int_{B_{R_0}}\zeta\rho_R[\mathrm{I_4}].
\end{align*}
Let $\xi(r)=r^m$, where $r=|x|$, $m$ is some positive integer
sufficiently large. Through elementary computations, we get
\begin{align*}
    &\mathrm{I_4}\\
        =&\frac12\frac{\xi'^2(r)}{r^2}\sum_{i,j=1}^n\left(GQG^T\right)_{ij}x_ix_j\\
        &-\frac12\left[m(m-2)r^{m-4}\sum_{i,j=1}^n\left(GQG^T\right)_{ij}x_ix_j\right.\\
&\phantom{-\frac12[}\left.+mr^{m-2}\textup{Tr}\left(GQG^T\right)\right]-mr^{m-2}(f-D_wK)\cdot x+N\\
        \geq&-\frac12\left[nm(m-2)\left|\left|GQG^T\right|\right|_\infty+m\,\textup{Tr}\left(GQG^T\right)\right]r^{m-2}\\
        &-mr^{m-2}(f-D_wK)\cdot x+N\geq C',
\end{align*}
where $C'$ is a positive constant independent of $R_0$, by condition
(\ref{condition c'}). For large enough $m$, we have
\begin{align*}
   |\mathrm{I_3}|
        &\leq e^{-R_0^m}\\
	&\cdot\int_{B_R}\left|n\left|\left|GQG^T\right|\right|_\infty\left[\frac12m^2r^{2m-2}-m\left(\frac12m-1\right)r^{m-2}\right]\right.\\
        &\phantom{\int_{B_R}|}\left.-\frac12m\,\textup{Tr}\left(GQG^T\right)r^{m-2}-m(f-D_wK)\cdot
        xr^{m-2}\right|\rho_R\\
        &\leq e^{-R_0^m}\left(\frac12nm(m+1)\left|\left|GQG^T\right|\right|_\infty
    R_0^{2m-2}+C''\right)\int_{B_R}\rho_R\\
        &\leq \left(\frac12nm(m+1)\left|\left|GQG^T\right|\right|_\infty R_0^{2m-2}+C''\right)\\
	&\phantom{\leq11}\cdot  e^{CT-R_0^m}\int_{B_R}\sigma_{0,R}\triangleq\gamma(R_0).
\end{align*}
The last inequality follows by the similar argument of (\ref{int
rho}). Hence,
\[
    \frac
    d{dt}\int_{B_{R_0}}\zeta\rho_R\geq-\gamma(R_0)-C'\int_{B_{R_0}}\zeta\rho_R.
\]
This implies
\begin{align}\label{zeta rho R}\notag
    &\int_{B_{R_0}}\zeta\rho_R(x,T)\\\notag
\geq& e^{-C'T}\int_{B_{R_0}}\zeta\sigma_{0,R}(x)+\frac{\gamma(R_0)}{C'}\left(e^{-C'T}-1\right)\\\notag
        \geq&e^{-C'T}\int_{B_{R_0}}\zeta\sigma_{0,R}(x)\\\notag
            &+\left(\frac12nm(m+1)\left|\left|GQG^T\right|\right|_\infty
    R_0^{2m-2}+C''\right)\\
	&\phantom{+a}\cdot\frac{e^{(C-C')T-R_0^m}}{C'}\left(1-e^{C'T}\right)\int_{B_R}\sigma_{0,R}(x).
\end{align}
Let $R_0\rightarrow\infty$, we have
\[
    \int_{\mathbb{R}^n}e^{-|x|^m}\rho(x,T)\geq
    e^{-C'T}\int_{\mathbb{R}^n}e^{-|x|^m}\sigma_0(x).
\]
\end{IEEEproof}

\section{Numerical simulations}

In this section, we shall apply our algorithm to both ``time-invariant" case and ``time-varying" case. The numerical simulations support our theorems. In our implementation, we adopt the Hermite spectral method (HSM) to get the approximate solution of \eqref{Kolmogorov equation}. Thus, the basis functions $\{\phi_l\}_{l=1}^\infty$ in \eqref{proj to basis} are 
choosen to be the generalized Hermite functions $\{H_n^{\alpha,\beta}(x)\}_{n=0}^\infty$. We refer the interested readers to the detailed definitions in \cite{LY}. 

For $N>0$, let us denote $\mathcal{R}_N$ the subspace spanned by the first $N$ generalized Hermite functions:
\begin{align*}
    \mathcal{R}_N=&\textup{span}\{H_0^{\alpha,\beta}(x),\cdots,H_N^{\alpha,\beta}(x)\}.
\end{align*}

The formulation of HSM to (\ref{Kolmogorov equation}) in 1-dimension is to find $u_N(x,t)\in\mathcal{R}_N$ such that
\begin{align}\label{weak formulation}
    \left\{ \begin{aligned}
       \langle\partial_tu_N(x,t),\varphi\rangle=&-\frac12\langle\partial_x[(GQG^T)u_N],\partial_x\varphi\rangle\\
	&+\langle fu_N,\partial_x\varphi\rangle-\frac12\langle(h^TS^{-1}h)u_N,\varphi\rangle\\     
       u_N(x,0)=&P_Nu_0(x),
    \end{aligned} \right.
\end{align}
for any $\varphi\in\mathcal{R}_N$, where $\langle\cdot,\cdot\rangle$ denotes the scalar product in $L^2(\mathbb{R})$ and $P_N$ is the projection operator such that $P_N:\,L^2(\mathbb{R})\rightarrow\mathcal{R}_N$. Write the solution $u_N\in\mathcal{R}_N$ in the form
\begin{align*}
    u_N(x,t)=\sum_{n=0}^Na_n(t)H_n^{\alpha,\beta}(x),
\end{align*}
and take the test function $\varphi\in\mathcal{R}_N$ in (\ref{weak formulation}) to be $H_n^{\alpha,\beta}$, $n=0,\cdots,N$. From (\ref{weak formulation}) and the properties of generalized Hermite functions, $\vec{a}(t):=(a_0(t), a_1(t),\cdots,a_{N}(t))^T$ satisfies the ODE
\begin{align}\label{ODE}
    \partial_t\vec{a}(t)=A\vec{a}(t),
\end{align}
where $A$ is a $(N+1)\times(N+1)$ matrix, may depend on $t$, if $G$, $Q$, $f$ or $h$ is explicitly time-dependent. This ODE can be precomputed. The only difference between ``time-varying" case and ``time-invariant" case is that it costs much more memory to store the off-line data in the ``time-varying" case, as we explained before in section II. Nevertheless, it doesn't change the off-line virture of our algorithm. We refer the interested readers of the implementation to \cite{LY}, and we shall omit the technical details in this paper.

Once $u_N$ at each step is obtained, $\hat{\rho}$ can be recovered by
(\ref{Rozovsky's reverse transformation}). The conditional expectation of the
state $x_t$ is computed by definition
\[
    \mathbb{E}\left[x,\{y_\tau\}_{0\leq\tau\leq t}\right](t) =
    \frac{\int_\mathbb{R}x\hat{\rho}(x,t)dx}{\int_\mathbb{R}\hat{\rho}(x,t)dx}.
\]

\subsection{``time-invariant" case: the 1D cubic sensor}

Let us consider the following model
\begin{equation*}
   \left\{ \begin{aligned}
        dx_t &= dv_t\\
        dy_t &= x_t^3dt+dw_t,
\end{aligned} \right.
\end{equation*}
where $x_t$, $y_t\in\mathbb{R}$, $v_t$, $w_t$ are scalar Brownian
motion processes with $E[dv_t^Tdv_t]=1$, $E[dw_t^Tdw_t]=1$. The 1D
Kolmogorov forward equation (\ref{Kolmogorov equation}) here is
\begin{align}\label{cubic sensor}
        u_t=\frac12u_{xx}-\frac12x^6u,
\end{align}
at each time step. We assume the inital density function $u_0(x)=e^{-x^4/4}$ and the updated initial data are
\[
    u_i(x,\tau_i)=e^{x^3\cdot dy_t }u_{i-1}(x,\tau_i).
\]
\begin{figure}[!t]
          \centering
	\includegraphics[trim = 30mm 75mm 30mm 70mm, clip, scale=0.5]{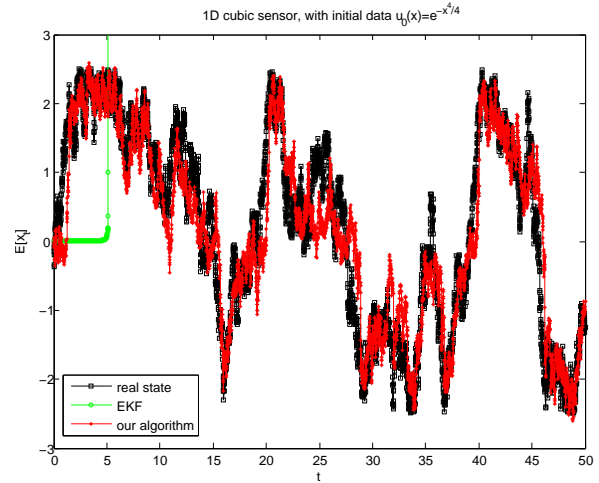}
          \caption{1D cubic sensor, with the initial condition $u_0(x)=e^{-\frac{x^4}4}$. Black: real state; Green: extended Kalman filter; Red: our algorithm.}
          \label{fig-cubic}
\end{figure}
In Figure \ref{fig-cubic}, we see that our algorithm tracks the state's expectation very well, while the
extended Kalman filter (EKF) completely fails around $t=5$. The total simulation time is $T=50$, and the update time step is $dt=\tau_{i+1}-\tau_i=0.01$. It costs our algorithm only around $4.88s$ to finish the simulation, i.e. the updated time is less than $10^{-3}s$.

\subsection{``time-varying" case: the 1D almost linear sensor}

\begin{figure}[!t]
          \centering
	\includegraphics[trim = 30mm 85mm 30mm 80mm, clip, scale=0.55]{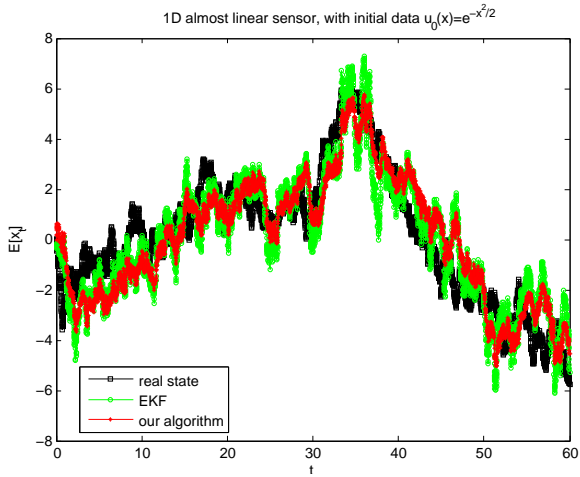}
          \caption{1D ``time-varying" almost linear sensor, with the initial condition $u_0(x)=e^{-\frac{x^2}2}$. Black: real state; Green: extended Kalman filter; Red: our algorithm. }
          \label{fig-coslineartime}
\end{figure}

The 1D almost linear sensor we are considering is 
\begin{align}
	\left\{\begin{aligned}
		dx_t =& [1+0.1\cos{(20\pi t)}]dv_t\\
		dy_t =&x_t[1+0.25\cos{(x_t)}]dt+dw_t,
	\end{aligned}\right.
\end{align}
where $x_t$, $y_t\in\mathbb{R}$, $v_t$, $w_t$ are scalar Brownian motion processes with $E[dv_t^Tdv_t]=E[dw_t^Tdw_t]=1$. The Kolmogorov forward equation \eqref{Kolmogorov equation} in this example is 
\begin{align*}
	u_t = \frac12[1+0.1\cos{(20\pi t)}]^2u_{xx}-\frac12x^2[1+0.25\cos{(x)}]^2u,
\end{align*}
with the initial data $u_0(x)=e^{-x^2/2}$ and the updated initial data
\begin{align*}
	u_i(x,\tau_i)=e^{x^2[1+0.25\cos{(x)}]\cdot dy_t}u_{i-1}(x,\tau_i),
\end{align*}
$i=1,2,\cdots,k$. In Figure \ref{fig-coslineartime}, our algorithm tracks the state's expectation at least as well as the EKF. The total simulation time is $T=60$, and the update time step is $dt=\tau_{i+1}-\tau_i=0.01$. It costs our algorithm only around $3.17s$ to complete the simulatoin, i.e. the updated time is less than $5\times10^{-4}s$.

\section{Conclusion}

In this paper, we extend the algorithm developed in \cite{YY} to the most general nonlinear filterings. We theoretically verified
that under very mild growth conditions on the drift term and the observation term, the unique non-negative weak
solution $\rho$ of its associated ``pathwise-robust" DMZ equation can be approximated by the
solution $\rho_R$ of the DMZ equation restricted on a large
ball $B_R$ with 0-Dirichlet boundary condition. The error
of this approximation tends to zero exponentially as the
radius of the ball $R$ approaching infinity. Moreover, $\rho_R$
can be efficiently approximated by our algorithm. We
show that the approximate solution $\hat{\rho}_R$ obtained by our
algorithm converges to $\rho_R$ in the $L^1$ sense for all $t\in[0,T]$, as the partition of time becomes finer, and a precise error estimate
of this convergence is given explicitly. Equally important, our algorithm preserves the two advantages of KF: ``memoryless" and ``real time". We also give the detail explanation of the off-line virture of our algorithm in the formulation. Numerical experiments support the feasibility and efficiency of our
algorithm.


%

\appendices
\section{Existence and uniqueness of the solution}

\renewcommand{\thetheorem}{A.\arabic{theorem}}
\renewcommand{\theequation}{A.\arabic{equation}}
\renewcommand{\theremark}{A.\arabic{remark}}
\setcounter{theorem}{0}
\setcounter{remark}{0}

\setcounter{equation}{0}
Before we show the existence of the weak solution, we shall give
a priori estimations of up to the first order derivative of the
solution to the robust DMZ equation on $B_R\times[0,T]$.
\begin{theorem}\label{a priori estimate}
    Consider the ``pathwise-robust" DMZ equation (\ref{robust DMZ eqn on ball})
    on $\mathbb{Q}_R:=B_R\times[0,T]$, where $B_R=\{x\in\mathbb{R}^n:|x|\leq R\}$ is a ball of radius
    $R$. Assume that
    \begin{align}\label{d/dt (GQG^T)}
        \left\|\frac d{dt}(GQG^T)\right\|_\infty<\infty,
    \end{align}
    for all $t\in[0,T]$. Suppose there exists a positive function $g(x)$ on
    $\mathbb{R}^n$ such that for all $t\in[0,T]$, $g$ and $\tilde{g}\triangleq g+\log{|D_wJ|}$ satisfy
\begin{enumerate}
     \item $\left|D_wg+\frac12\nabla(GQG^T)-F\right|^2+2\lambda_1J\leq
            C,$\inlineeq{}\label{A condition 1}
     \item $D_w^2g+2D_wg\cdot\nabla
        g+2[\nabla(GQG^T)-F]\cdot\nabla
        g+\frac12\nabla^2(GQG^T)-\textup{div}F+J\leq C,$\inlineeq{}\label{A condition 2}
     \item $D_w^2\tilde{g}+2D_w\tilde{g}\cdot\nabla
        \tilde{g}+2[\nabla(GQG^T)-F]\cdot\nabla
        \tilde{g}+\frac12\nabla^2(GQG^T)-\textup{div}F+J\leq
        C,$\inlineeq{}\label{A condition 3}
     \item $\int_{\mathbb{R}^n}e^{2\tilde{g}}\sigma^2(x)\leq C\quad\textup{and}\quad\int_{\mathbb{R}^n}e^{2g}D_w\sigma\cdot\nabla\sigma\leq
     C,$\inlineeq{}\label{A condition 4}
\end{enumerate}
     where $C$ is a generic constant, which may differ from line to
     line, and $\nabla(*)=\left[\sum_{i=1}^n\frac{\partial (*)_{ij}}{\partial
     x_i}\right]_{j=1}^n$, $\nabla^2(*)=\sum_{i,j=1}^n\frac{\partial^2 (*)_{ij}}{\partial x_i\partial
     x_j}$. Then, for $0\leq t\leq T$,
    \begin{align}\label{a priori of order zero}
        \int_{B_R}e^{2g}\rho_R^2(x,t)dx\leq e^{Ct}\int_{B_R}e^{2g}\sigma^2(x)dx,
\end{align}
\begin{align}\label{a priori of order
            one}\notag
        &\int_{B_R}e^{2g}D_w\rho_R(x,t)\cdot\nabla\rho_R(x,t)dx\\
            \leq&e^{Ct}\int_{B_R}e^{2g}D_w\sigma(x)\cdot\nabla\sigma(x)dx+Ce^{Ct}\int_{B_R}e^{2\tilde{g}}\sigma^2(x)dx,
    \end{align}
    where $D_w$ and $J(x,t)$ is defined in (\ref{Dw})and (\ref{J}), respectively.
\end{theorem}
\begin{remark}
    The conditions in Theorem {\ref{a priori estimate}} are easily checked, if the drift terms $h(x)$ and $f(x)$ are at most polynomial
    growth in $r=|x|$. However, in general, the existence of such
    $g$ is not always available.
\end{remark}

\noindent{\bf Proof:}\quad Let $g$ be some positive function on
$\mathbb{R}^n$.
\begin{align}\label{A 1}\notag
    \frac d{dt}\int_{B_R}e^{2g}\rho_R^2
        =&\int_{B_R}e^{2g}\rho_RD_w^2\rho_R+2\int_{B_R}e^{2g}\rho_R(F\cdot\rho_R)\\\notag
           & +2\int_{B_R}e^{2g}J\rho_R^2\\
        \triangleq&\mathrm{I}+\mathrm{II}+\mathrm{III}.
\end{align}
Apply integration by parts to $\mathrm{I}$ and $\mathrm{II}$ in
(\ref{A 1})
\begin{align*}
    \mathrm{I}
        =&-2\int_{B_R}\rho_Re^{2g}D_wg\cdot\nabla\rho_R
            -\int_{B_R}e^{2g}D_w\rho_R\cdot\nabla\rho_R\\
           & -\int_{B_R}e^{2g}\rho_R\nabla(GQG^T)\cdot\nabla\rho_R\\
        \leq&-2\int_{B_R}\rho_Re^{2g}D_wg\cdot\nabla\rho_R
            -\int_{B_R}e^{2g}\rho_R\nabla(GQG^T)\cdot\nabla\rho_R\\
        &\triangleq\mathrm{I_1}+\mathrm{I_2}.
\end{align*}
Integration by parts further, we have
\begin{align}\label{A 2}\notag
    \mathrm{I_1}
        =&4\int_{B_R}e^{2g}\rho_R^2D_wg\cdot\nabla g
            +2\int_{B_R}e^{2g}\rho_RD_w g\cdot\nabla\rho_R\\
           &+2\int_{B_R}e^{2g}\rho_R^2\nabla(GQG^T)\cdot\nabla g
         +2\int_{B_R}e^{2g}\rho_R^2D_w^2g.
\end{align}
Notice that the second term of the right-hand side of (\ref{A 2}) is
$-\mathrm{I_1}$, we have
\begin{align}\label{I_1}\notag
    \mathrm{I_1}
        =&2\int_{B_R}e^{2g}\rho_R^2D_wg\cdot\nabla g\\
            &+\int_{B_R}e^{2g}\rho_R^2[\nabla(GQG^T)\cdot\nabla
            g+D_w^2g].
\end{align}
The similar argument applies to $\mathrm{I_2}$:
\begin{align}\label{I_2}\notag
    \mathrm{I_2}
        =&\int_{B_R}e^{2g}\rho_R^2\nabla(GQG^T)\cdot\nabla g\\
            &+\frac12\int_{B_R}e^{2g}\rho^2_R\nabla^2(GQG^T).
\end{align}
Thus,
\begin{align}\label{I}
    \mathrm{I}
	\leq&\int_{B_R}e^{2g}\rho_R^2\left[D_w^2g+2D_wg\cdot\nabla g\right.
\\
&\phantom{\int_{B_R}e^{2g}\rho_R^2[a}\left.+2\nabla(GQG^T)\cdot\nabla
    g+\frac12\nabla^2(GQG^T)\right].
\end{align}
The same trick of $\mathrm{I_1}$ applies to $\mathrm{II}$ in (\ref{A
1}), we obtain
\begin{align}\label{II}
    \mathrm{II}
        =-\int_{B_R}e^{2g}\rho_R^2[2F\cdot\nabla g+\textup{div}F].
\end{align}
Substitute (\ref{I}) and (\ref{II}) back to (\ref{A 1}), we obtain
\begin{align*}
    \frac d{dt}&\int_{B_R}e^{2g}\rho_R^2\\
        \leq&\int_{B_R}e^{2g}\rho_R^2\left\{D_w^2g
            +2D_wg\cdot\nabla g+2[\nabla(GQG^T)-F]\cdot\nabla g\right.\\
            &\phantom{\int_{B_R}e^{2g}\rho_R^2[}\left.+\frac12\nabla^2(GQG^T)-\textup{div}F+J\right\}\\
        \leq&C\int_{B_R}e^{2g}\rho_R^2,
\end{align*}
by condition (\ref{A condition 2}). (\ref{a priori of order zero})
follows directly from Gronwall's inequality. To show (\ref{a priori
of order one}), we consider
\begin{align}\label{A 6}\notag
    \frac d{dt}&\int_{B_R}e^{2g}D_w\rho_R\cdot\nabla\rho_R\\\notag
        =&\int_{B_R}e^{2g}\sum_{i,j=1}^n\frac d{dt}(GQG^T)_{ij}\frac{\partial\rho_R}{\partial x_i}\frac{\partial\rho_R}{\partial x_j}\\\notag
            &+2\int_{B_R}e^{2g}\sum_{i,j=1}^n(GQG^T)_{ij}\frac{\partial}{\partial x_i}\left(\frac{\partial\rho_R}{\partial
                t}\right)\frac{\partial\rho_R}{\partial x_j}\\
        \triangleq&\mathrm{IV}+\mathrm{V}.
\end{align}
Due to condition (\ref{d/dt (GQG^T)}), $\mathrm{IV}$ of (\ref{A 6})
turns out to be
\begin{align}\label{IV}\notag
    \mathrm{IV}
        \leq&\frac12\left\|\frac d{dt}(GQG^T)\right\|_{\infty}\int_{B_R}e^{2g}\sum_{i,j=1}^n\left[\left(\frac{\partial\rho_R}{\partial
            x_i}\right)^2+\left(\frac{\partial\rho_R}{\partial
            x_j}\right)^2\right]\\\notag
        &=n\left\|\frac
        d{dt}(GQG^T)\right\|_{\infty}\int_{B_R}e^{2g}|\nabla\rho_R|^2\\
        \leq&\frac n{\lambda_1}\left\|\frac d{dt}(GQG^T)\right\|_{\infty}\int_{B_R}e^{2g}D_w\rho_R\cdot\nabla\rho_R,
\end{align}
since $D_w\rho_R\cdot\nabla\rho_R\geq\lambda_1|\nabla\rho_R|^2$.
Next, $\mathrm{V}$ in (\ref{A 6}) is
\begin{align}\label{A 7}\notag
    \mathrm{V}
        =&-2\int_{B_R}e^{2g}[(2D_wg+\nabla(GQG^T))\cdot\nabla\rho_R+D_w^2\rho_R]\\\notag
&\phantom{-2\int_{B_R}e^{2g}[}\cdot\left(\frac12D_w^2\rho_R+F\cdot\nabla\rho_R+J\rho_R\right)\\\notag
        =&-\int_{B_R}e^{2g}\left\{D_w^2\rho_R+\left[D_wg+\frac12\nabla(GQG^T)+F\right]\cdot\nabla\rho_R\right\}^2\\\notag
         &+\int_{B_R}e^{2g}\left[D_wg+\frac12\nabla(GQG^T)-F\right]^2|\nabla\rho_R|^2\\\notag
         &-2\int_{B_R}e^{2g}[D_w^2\rho_R+(2D_wg+\nabla(GQG^T))\cdot\nabla\rho_R]J\rho_R\\\notag
        \leq&\int_{B_R}e^{2g}\left[D_wg+\frac12\nabla(GQG^T)-F\right]^2|\nabla\rho_R|^2\\
         &-2\int_{B_R}e^{2g}[D_w^2\rho_R+(2D_wg+\nabla(GQG^T))\cdot\nabla\rho_R]J\rho_R.
\end{align}
Notice that
\begin{align}\label{A 8}\notag
    &\int_{B_R}e^{2g}D_w^2\rho_RJ\rho_R\\\notag
        =&-\int_{B_R}e^{2g}\left[2(D_wg\cdot\nabla\rho_R)J\rho_R+JD_w\rho_R\cdot\nabla\rho_R\right.\\
        &\phantom{-\int_{B_R}e^{2g}[}\left.+(D_w\rho_R\cdot\nabla
        J)\rho_R+\nabla(GQG^T)\cdot\nabla\rho_RJ\rho_R\right].
\end{align}
Take (\ref{A 8}) into account, $\mathrm{V}$ becomes
\begin{align}\label{A 9}\notag
    \mathrm{V}\leq&\int_{B_R}e^{2g}\\\notag
        &\cdot\left\{\frac1{\lambda_1}\left\{\left[D_wg+\frac12\nabla(GQG^T)-F\right]^2+1\right\}
               +2J\right\}\\
	&\cdot D_w\rho_R\cdot\nabla\rho_R+\int_{B_R}e^{2g}|D_wJ|^2\rho_R^2.
\end{align}
Combine (\ref{IV}) and (\ref{A 9}), we have
\begin{align}\notag
    \frac d{dt}&\int_{B_R}e^{2g}D_w\rho_R\cdot\nabla\rho_R\\\notag
        \leq&\int_{B_R}e^{2g}\left\{\frac1{\lambda_1}\left\{n\left\|\frac d{dt}(GQG^T)\right\|_\infty\right.\right.\\\notag
		&\phantom{\int_{B_R}e^{2g}}\left.\left.
                +\left[D_wg+\frac12\nabla(GQG^T)-F\right]^2+1\right\}
                +2J\right\}\\\notag
            &\phantom{\int_{B_R}e^{2g}}\cdot D_w\rho_R\cdot\nabla\rho_R\\
	&+\int_{B_R}e^{2g}|D_wJ|^2\rho_R^2.
\end{align}
By conditions (\ref{A condition 1})-(\ref{A condition 4}), the
estimate (\ref{a priori of order one}) follows
immediately.\hfill{$\Box$}

\noindent{\bf Proof of existence in Theorem \ref{existence and uniqueness}:}\quad
Let $R_k$ be a sequence of positive number such that
$\lim_{k\rightarrow\infty}R_k=\infty$. Let $\rho_k(x,t)$ be the
solution of the ``pathwise-robust" DMZ equation (\ref{robust DMZ eqn on ball})
on $B_{R_k}\times[0,T]$, where $B_{R_k}=\{x\in\mathbb{R}^n:|x|\leq
R_k\}$ is a ball of radius $R_k$. In view of Theorem \ref{a priori
estimate}, the sequence $\{\rho_k\}$ is a bounded set in
$H_0^{1;1}(\mathbb{Q}_{R_k})$. Thus, there exists a subsequence
$\{\rho_{k'}\}$ which is weakly convergent to $\rho$. Moreover,
$\rho$ has the weak derivative $\frac{\partial\rho}{\partial x_i}\in
L^2(\mathbb{Q}_{R_k})$, and $\frac{\partial\rho_{k'}}{\partial x_i}$
weakly tends to it. Now we claim that the weak derivative
$\frac{\partial\rho}{\partial t}$ exists. To see this, let
$\Phi(x,t)\in H_0^{1;1}(\mathbb{Q}_{R_k})$, then
\begin{align*}
    &\iint_{\mathbb{Q}_{R_k}}\frac12\sum_{i,j=1}^n(GQG^T)_{ij}\frac{\partial\Phi}{\partial
        x_j}\frac{\partial\rho}{\partial x_i}\\
        &\phantom{\iint_{\mathbb{Q}_{R_k}}}+\left[\sum_{i=1}^n\left(\sum_{j=1}^n\frac{\partial(GQG^T)_{ij}}{\partial x_j}-F_i\right)\frac{\partial\rho}{\partial
        x_i}-J\rho\right]\Phi\\
       =&\lim_{k'\rightarrow\infty}\iint_{\mathbb{Q}_{R_k}}\frac12\sum_{i,j=1}^n(GQG^T)_{ij}\frac{\partial\Phi}{\partial
        x_j}\frac{\partial\rho_{k'}}{\partial x_i}\\
        &+\left[\sum_{i=1}^n\left(\sum_{j=1}^n\frac{\partial(GQG^T)_{ij}}{\partial x_j}-F_i\right)\frac{\partial\rho_{k'}}{\partial
        x_i}-J\rho_{k'}\right]\Phi\\
       =&-\lim_{k'\rightarrow\infty}\iint_{\mathbb{Q}_{R_k}}\frac{\partial\rho_{k'}}{\partial
        t}\Phi
       =\lim_{k'\rightarrow\infty}\iint_{\mathbb{Q}_{R_k}}\rho_{k'}\frac{\partial\Phi}{\partial
        t}\\
       =&\iint_{\mathbb{Q}_{R_k}}\rho\frac{\partial\Phi}{\partial t}.
\end{align*}
Clearly,
$\rho(x,0)=\lim_{k'\rightarrow\infty}\rho_{k'}(x,0)=\sigma_0(x)$.

\begin{theorem}\label{uniqueness}
 	Assume further that for some $c>0$,
    \begin{align}\label{unique-cond1}
        \sup_{0\leq t\leq
        T}\int_{\mathbb{R}^n}e^{cr}\rho^2(x,t)dx<\infty,
    \end{align}
and
    \begin{align}\label{unique-cond2}
        \int_{\mathbb{Q}_T}|\nabla\rho(x,t)|^2dxdt<\infty,
    \end{align}
where $r=|x|$. Suppose that there exists a finite number $\alpha>0$
such that
    \begin{align}\label{alpha}
        2J(x,t)-\frac1{4\lambda_1}[cD_wr-(F(x,t)+\tilde{F}(x,t))]^2\leq\alpha,
    \end{align}
for all $(x,t)\in\mathbb{Q}_T$, where $\lambda_1$ is the smallest
eigenvalue of the matrix $(GQG^T)$,
    \begin{align}\label{tilde F}
        \tilde{F}(x,t)=\left[\frac12\sum_{j=1}^n(GQG^T)_{ij}+\sum_{j=1}^n(GQG^T)_{ij}\frac{\partial K}{\partial
        x_j}-f_i\right]_{i=1}^n,
    \end{align}
and $J(x,t)$ is defined as in (\ref{J}). Then the non-negative weak solution $\rho(x,t)$ of the ``pathwise-robust" DMZ equation on $\mathbb{Q}_T$ is unique.
\end{theorem}

\noindent{\bf Proof of uniqueness of Theorem \ref{existence and uniqueness} (Theorem \ref{uniqueness}):} To show the uniqueness of the solution, we only need to show that
$\rho(x,t)=0$ on $\mathbb{Q}_T$ if $\rho(x,0)=0$. Let $\alpha T<1$.
For any test function $\psi(x,t)=e^{cr}\Phi(x,t)$, where $r=|x|$,
$c$ is some constant and $\Phi(x,t)\in H_0^{1;1}(\mathbb{Q}_T)$,
then $\rho(x,t)$ satisfies
\begin{align}\label{C 1}\notag
    &\int_{\mathbb{R}^n}\rho(x,T)\Phi(x,T)e^{cr}dx
         -\int_0^T\int_{\mathbb{R}^n}\rho(x,t)\frac{\partial\Phi}{\partial
         t}(x,t)e^{cr}dxdt\\\notag
    =&\int_{\mathbb{Q}_T}-\frac12e^{cr}\nabla\Phi(x,t)\cdot D_w\rho(x,t)
      -\frac c2e^{cr}\Phi(x,t)\nabla r\cdot D_w\rho(x,t)\\\notag
      &\phantom{\int_{\mathbb{Q}_T}}+\tilde{F}(x,t)\cdot\nabla\rho(x,t)\Phi(x,t)e^{cr}\\
      &\phantom{\int_{\mathbb{Q}_T}}+J(x,t)\rho(x,t)\Phi(x,t)e^{cr}dxdt.
\end{align}
where $\tilde{F}$ is defined in (\ref{tilde F}). Approximate
$\rho(x,t)$ by $\Phi(x,t)$ in the $H^{1;1}(\mathbb{Q}_T)$-norm, we
get
\begin{align}\label{C 4}\notag
    &\int_{\mathbb{R}^n}\rho^2(x,T)e^{cr}dx\\\notag
    =&\int_{\mathbb{Q}_T}e^{cr}\left[-D_w\rho(x,t)\cdot\nabla\rho(x,t)
            -c\rho(x,t)\nabla r\cdot D_w\rho(x,t)\right.\\\notag
         &\phantom{\int_{\mathbb{Q}_T}e^{cr}[a}+(\tilde{F}(x,t)+F(x,t))\cdot\nabla\rho(x,t)\rho(x,t)\\\notag
           & \phantom{\int_{\mathbb{Q}_T}e^{cr}[a}\left.+2J(x,t)\rho^2(x,t)\right]dxdt.\\\notag
         \leq&\int_{\mathbb{Q}_T}e^{cr}[-\lambda_1|\nabla\rho(x,t)|^2-c\rho(x,t)D_w
        r\cdot\nabla\rho(x,t)\\\notag
            &\phantom{\int_{\mathbb{Q}_T}e^{cr}[a}+(F(x,t)+\tilde{F}(x,t))\cdot\nabla\rho(x,t)\rho(x,t)\\\notag
                 &\phantom{\int_{\mathbb{Q}_T}e^{cr}[a}+2J(x,t)\rho^2(x,t)]dxdt.\\\notag
        =&-\lambda_1\int_{\mathbb{Q}_T}e^{cr}\left\{\frac1{2\lambda_1}[cD_wr-(F(x,t)+\tilde{F}(x,t))]\rho(x,t)\right.\\\notag
            &\phantom{-\lambda_1\int_{\mathbb{Q}_T}e^{cr}[aa}+\left.|\nabla\rho(x,t)|\right\}^2dxdt\\\notag
         &+\int_{\mathbb{Q}_T}e^{cr}\left\{2J(x,t)-\frac1{4\lambda_1}[cD_wr-(F(x,t)+\tilde{F}(x,t))]^2\right\}\\\notag
	&\phantom{\int_{\mathbb{Q}_T}e^{cr}}\cdot\rho^2(x,t)dxdt\\\notag
        \leq&\int_{\mathbb{Q}_T}e^{cr}\left\{2J(x,t)-\frac1{4\lambda_1}[cD_wr-(F(x,t)+\tilde{F}(x,t))]^2\right\}\\
	&\phantom{\int_{\mathbb{Q}_T}e^{cr}}\cdot\rho^2(x,t)dxdt,
\end{align}
due to the positive definite of $(GQG^T)$. By condition
(\ref{alpha}), we have
\begin{align}\label{C 6}
    \int_{\mathbb{R}^n}e^{cr}\rho^2(x,T)dx\leq\alpha\int_{\mathbb{Q}_T}e^{cr}\rho^2(x,t)dxdt.
\end{align}
According to the mean value theorem, there exists $T_1\in(0,T)$ such
that
\begin{align}\label{C 7}\notag
    \int_{\mathbb{Q}_T}e^{cr}\rho^2(x,t)dxdt
        =&\int_0^T\int_{\mathbb{R}^n}e^{cr}\rho^2(x,t)dxdt\\
        =&T\int_{\mathbb{R}^n}e^{cr}\rho^2(x,T_1)dx.
\end{align}
Apply (\ref{C 6}) and (\ref{C 7}) recursively, there exists
$T_m\in(0,T)$ such that
\[
    \int_{\mathbb{R}^n}e^{cr}\rho^2(x,T)dx
        \leq(\alpha T)^m\int_{\mathbb{R}^n}e^{cr}\rho^2(x,T_m)dx.
\]
Since $\alpha T<1$, we conclude that $\rho(x,t)\equiv0$ for a.e
$(x,t)\in\mathbb{Q}_T$.\hfill{$\Box$}

\section{Proof of Theorem \ref{coro}}

\renewcommand{\thetheorem}{B.\arabic{theorem}}
\renewcommand{\theequation}{B.\arabic{equation}}
\renewcommand{\theremark}{B.\arabic{remark}}
\setcounter{theorem}{0}
\setcounter{remark}{0}

\setcounter{equation}{0}

\begin{IEEEproof}[Proof of  Theorem \ref{coro}] Let $v=\rho-\rho_R$ as in the proof of
Theorem \ref{convergence of rho R to rho}. By the maximum principle,
we have that $v\geq0$ for all $(x,t)\in B_R\times[0,T]$. Choose the
test function $\psi$ in Lemma \ref{IBP} as
\[
    \Phi(x)=\gamma(x)\varrho(x),
\]
where $\gamma(x)=e^{\frac12\phi_1(x)}$ and $\phi_1(x)$, $\varrho(x)$
are defined in the proof of Proposition \ref{rho concentrate at
origin} and Theorem \ref{convergence of rho R to rho}. It follows
directly that $\Phi|_{\partial B_R}=\nabla_x\Phi|_{\partial B_R}=0$,
by the fact that $\varrho|_{\partial B_R}=\nabla\varrho|_{\partial
B_R}=0$. Apply Lemma \ref{IBP} to $v$ taking place of $\rho_\Omega$ with
the test function $\Phi$, we have
\begin{align*}
    &\frac{d}{dt}\int_{B_R}\Phi v\\
        =&\frac12\int_{B_R}D_w^2\Phi v+\int_{B_R}(f-D_wK)\cdot\Phi
            v+\int_{B_R}\Phi Nv\\
        =&\frac12\int_{B_R}(D_w^2\gamma\varrho+2D_w\gamma\cdot\nabla\varrho+\gamma
            D_w^2\varrho)v\\
	&+\int_{B_R}(f-D_wK)\cdot(\nabla\gamma\varrho+\gamma\nabla\varrho)v+\int_{B_R}\gamma\varrho Nv.
\end{align*}
All the boundary integrals vanish due to the similar arguments in
Theorem \ref{convergence of rho R to rho}. Recall that
$\gamma(x)=e^{\frac12\phi_1(x)}$ and
$\varrho(x)=e^{-\phi_2(x)}-e^{-R}$.  Direct computations yield that
\begin{align*}
    &\frac{d}{dt}\int_{B_R}\Phi v\\
        =&\frac12\int_{B_R}\left[\frac12e^{\frac12\phi_1}\left(D_w^2\phi_1+\frac12D_w\phi_1\cdot\nabla\phi_1\right)\varrho\right.\\
           & \phantom{\frac12\int_{B_R}[}-e^{\frac12\phi_1}D_w\phi_1\cdot e^{-\phi_2}\nabla\phi_2\\
         &\left.\phantom{\frac12\int_{B_R}[}+\gamma
         e^{-\phi_2}\left(D_w\phi_2\cdot\nabla\phi_2-D_w^2\phi_2\right)\right]v\\
         &+\int_{B_R}(f-D_wK)\cdot\left(\frac12e^{\frac12\phi_1}\nabla\phi_1\varrho-\gamma
         e^{-\phi_2}\nabla\phi_2\right)v\\
	&+\int_{B_R}\gamma\varrho
         Nv\\
        =&\int_{B_R}\Phi v\left[\frac14\left(D_w^2\phi_1+\frac12D_w\phi_1\cdot\nabla\phi_1\right)-\frac12D_w\phi_1\cdot\nabla\phi_2\right.\\
           & \phantom{\int_{B_R}\Phi v[a}+\frac12\left(D_w\phi_2\cdot\nabla\phi_2-D_w^2\phi_2\right)\\
         &\phantom{\int_{B_R}\Phi v[a}\left.+(f-D_wK)\cdot(\frac12\nabla\phi_1-\nabla\phi_2)+N\right]\\
         &+e^{-R}\int_{B_R}\gamma v\left[-\frac12D_w\phi_1\cdot\nabla\phi_2\right.\\
	&\phantom{e^{-R}\int_{B_R}\gamma v[aa}+\frac12\left(D_w\phi_2\cdot\nabla\phi_2-D_w^2\phi_2\right)\\
            &\phantom{e^{-R}\int_{B_R}\gamma v[aa}\left.-(f-D_wK)\cdot\nabla\phi_2\right]\\
        \triangleq&\int_{B_R}\Phi
        v[\mathrm{VI}]+e^{-R}\int_{B_R}\gamma v[\mathrm{VII}],
\end{align*}
By the similar estimates (\ref{estimate of
D_w^2phi_1})-(\ref{estimate of (f-D_wK)nabla phi_1}), (\ref{estimate
of D_w^2phi_2})-(\ref{estimate of (f-D_wK)nabla phi_2}), we have
\begin{align*}
    \sup_{B_R}|\mathrm{VI}|&\leq
    17n\left|\left|GQG^T\right|\right|_\infty+5|f-D_wK|+N,\\
    \sup_{B_R}|\mathrm{VII}|&\leq
    16n\left|\left|GQG^T\right|\right|_\infty+4|f-D_wK|.
\end{align*}
Hence,
\begin{align*}
    \frac{d}{dt}\int_{B_R}\Phi v
        \leq& C\int_{B_R}\Phi v+e^{-R}\tilde{C}\int_{B_R}e^{\phi_1}v\\
        \leq& C\int_{B_R}\Phi v+e^{-R}\tilde{C}\int_{B_R}e^{\phi_1}\rho\\
        \leq& C\int_{B_R}\Phi
        v+\tilde{C}e^{-R+Ct}\int_{B_R}e^{\phi_1}\sigma_0(x)\\
            \leq& C\int_{B_R}\Phi v+\tilde{C}e^{-R+Ct}\int_{\mathbb{R}^n}e^{\phi_1}\sigma_0(x),
\end{align*}
by condition (\ref{condition c1}), (\ref{condition c3}) and
(\ref{estimation1 of thm concentration at origin}). By the similar
argument in the proof of Theorem \ref{convergence of rho R to rho},
where we get the estimate of $\int_{B_R}\varrho v$, we have
\begin{equation}\label{Thm 1.3 Eqn 1}
    \int_{B_R}\Phi v(x,T)
        \leq
        Ce^{-R}\int_{\mathbb{R}^n}e^{\sqrt{1+|x|^2}}\sigma_0(x),
\end{equation}
where $C$ is a generic constant, which depends on $T$. Recall that
$\varrho(x)=e^{-R[-(|x|^2/R^2-1)^2+1]}-e^{-R}$, it implies that
\begin{equation}\label{Thm 1.3 Eqn 2}
    \int_{B_R}\Phi v(x,T)
        \geq \frac12 e^{-\frac7{16}R}\int_{B_{\frac R2}}\gamma
        v(x,T).
\end{equation}
Combine (\ref{Thm 1.3 Eqn 1}) and (\ref{Thm 1.3 Eqn 2}), we obtain
that
\[
    \int_{B_{\frac R2}}\gamma
    v(x,T)\leq Ce^{-\frac9{16}R}\int_{\mathbb{R}^n}e^{\sqrt{1+|x|^2}}\sigma_0(x).
\]
This implies that
\begin{align*}
    \int_{B_{\frac R2}}\gamma\rho(x,T)
        \leq&\int_{B_{\frac
        R2}}\gamma\rho_R(x,T)\\
	&+Ce^{-\frac9{16}R}\int_{\mathbb{R}^n}e^{\sqrt{1+|x|^2}}\sigma_0(x)\\
        \leq& C(1+e^{-\frac9{16}R})\int_{\mathbb{R}^n}e^{\sqrt{1+|x|^2}}\sigma_0(x),
\end{align*}
by (\ref{estimation1 of thm concentration at
    origin}). Let $R\rightarrow\infty$,
    \[
        \int_{\mathbb{R}^n}\gamma\rho(x,T)
            \leq C\int_{\mathbb{R}^n}e^{\sqrt{1+|x|^2}}\sigma_0(x).
    \]
Consider the integration outside the large ball $B_R$,
    \begin{align*}
        e^{\frac12\sqrt{1+R^2}}\int_{|x|\geq R}\rho(x,T)
            \leq& \int_{|x|\geq R}\gamma\rho(x,T)\\
	\leq& C\int_{\mathbb{R}^n}e^{\sqrt{1+|x|^2}}\sigma_0(x).
    \end{align*}
    Therefore, we reach the conclusion that
    \[
        \int_{|x|\geq R}\rho(x,T)\leq
        Ce^{-\frac12\sqrt{1+|R|^2}}\int_{\mathbb{R}^n}e^{\sqrt{1+|x|^2}}\sigma_0(x).
    \]
\end{IEEEproof}



\ifCLASSOPTIONcaptionsoff
  \newpage
\fi

\end{document}